\documentclass[a4paper]{easychair}

\usepackage{graphicx}
\usepackage{fixltx2e}                    
\usepackage[utf8]{inputenc}            
\usepackage{amssymb,latexsym, dsfont, amsfonts, amsbsy, amsmath, mathrsfs, dsfont}            
\usepackage{mathtools}                   
\usepackage{bm}                          
\usepackage{enumerate}                   
\usepackage{verbatim}                    
\usepackage{url}                         
\usepackage{tabularx}
\usepackage{microtype}
\usepackage{latexsym}
\usepackage{amsmath,amssymb,amscd}
\usepackage{graphicx}
\usepackage{tikz}
\usepackage{enumitem}

\usepackage[all,cmtip]{xy}

\newtheorem{theorem}{Theorem}[section]
\newtheorem{corollary}[theorem]{Corollary}
\newtheorem{lemma}[theorem]{Lemma}
\newtheorem{proposition}[theorem]{Proposition}
\newtheorem{definition}[theorem]{Definition}
\newtheorem{remark}[theorem]{Remark}

\newcommand{\Go}[1]{{G\"{o}del} }

\usepackage{color}

\newcommand{\termDef}[1]{{\textbf{\textit{#1}}}}
\newcommand{\m}[1]{\mathfrak{#1}}
\DeclareMathSymbol{\sm}{\mathbin}{AMSa}{"39}

\renewcommand{\to}{\mathop{\rightarrow}} 

\newcommand{\alg}[1]{\mathbf{#1}}

\begin{document}
\setcounter{page}{1}     




\title{Axiomatization of crisp G\"odel modal logic}
\author{
	Ricardo O. Rodriguez$^{1,2}$  \and Amanda Vidal$^3$
}

\institute{{UBA. FCEyN. Departamento de Computaci\'on.} \\
\and
{CONICET-UBA. Inst. de Invest. en Cs. de la Computaci\'on.}\\
 Buenos Aires, Argentina.
\email{ricardo@dc.uba.ar}
	\and 
	Institute of Computer Sciences of the Czech Academy of Sciences, 
	Czech Republic\\
	\email{amanda@cs.cas.cz}
}

\authorrunning{R. Rodriguez, A. Vidal}
\titlerunning{Axiomatization of crisp G\"odel modal logic}

	\maketitle




\begin{abstract}
In this paper we consider the modal logic with both $\Box$ and $\Diamond$ arising from Kripke models with a crisp accessibility and whose propositions are valued over the standard G\"odel algebra $[0,1]_G$. We provide an axiomatic system extending the one from \cite{CaRo15} for models with a valued accessibility with Dunn axiom from positive modal logics, and show it is strongly complete with respect to the intended semantics. The axiomatizations of the most usual frame restrictions are given too. We also prove that in the studied logic it is not possible to get $\Diamond$ as an abbreviation of $\Box$, nor vice-versa, showing that indeed the axiomatic system we present does not coincide with any of the mono-modal fragments previously axiomatized in the literature.
\end{abstract}

\section{Introduction}\label{intro}
G{\"o}del Kripke models ({\bf GK}-models, or $\m{GK}$) are the generalization of the classical Kripke semantics for modal logics 
where both, the propositions at each world and the accessibility relation, are valued in the standard G{\"o}del algebra $[0,1]$.
A particular subclass of G{\"o}del Kripke models is the one given by the crisp ones ($\textbf{GK}^c$-models, or $\m{GK}^c$) where the accessibility relation only takes classical values, i.e. in $\{0,1\}$. This is to say, the frames are classical and is the models which incorporate some many-valued characteristics.

More general approaches, focusing mainly on finite residuated lattices, have been developed by Fitting~\cite{Fi92a, Fi92b}, Priest~\cite{Pr08b}, and Bou et al.~\cite{BoEsGoRo11}, and for other main fuzzy logics in \cite{HaTe13} and \cite{ViEsGo16}.
We remark that this approach to modal many-valued logics, starting from a Kripke semantics that behaves with respect to the FO semantics of the corresponding many-valued logic in the analogous way to how it does in the classical case (i.e., modalities can be translated to restricted quantifiers), differs from another main framework for modal substructural logics studied for instance in \cite{O93}, \cite{Re93b}, \cite{Kam02}.
This second approach, in contrast, can be seen as arising from a syntactic definition  of the logics, given by considering extensions of substructural logics with modalities governed by some of the usual axioms/rules of the modalities from classical modal logic. These logics enjoy completeness with respect to Kripke semantics different from the ones studied in this paper (which have additional relations, instead of considering valuated worlds), more similar to the ones appearing in modal intuitionistic logics.

The minimum logics over {\bf GK}-models have been investigated in some detail by Caicedo and Rodriguez~\cite{CaRo10,CaRo15} and Metcalfe and Olivetti~\cite{MeOl09,MeOl11}. Axiomatizations $\mathcal{GK}_\Box$ and $\mathcal{GK}_\Diamond$ were proposed for the logics with only one modality (box and diamond fragments) arising from $\m{GK}$ in~\cite{CaRo10}. There, it was also proved that the box fragment is not able to discriminate between crisp and non-crisp models, i.e.
the box fragment of {\bf GK}-models coincides with the box fragment of $\textbf{GK}^c$-models. On the contrary, the sets of valid
diamond formulae under both semantics are different. Interestingly enough, in \cite{CaRo10} it is proven that the diamond fragment
(over the whole class of models) enjoys the finite model property (FMP) with respect to its Kripke semantics, and its decidability
is established. This finite model property fails for both, the box fragment (for which, as we said above, both crisp and
general-valued Kripke semantics coincide)  and the diamond fragment of
$\textbf{GK}^c$-models. Nevertheless, decidability and PSPACE-completeness of validity in these $\Diamond$-fragments 
is established in \cite{MeOl09,MeOl11} using analytic Gentzen-style proof systems.

Further, in \cite{MeOl11}, the diamond fragment of $\textbf{GK}^c$-models is axiomatized ($\mathcal{GK}^c_\Diamond$) too.
It is subsequently shown in~\cite{CaRo15} that  the full logic with the two modalities arising from $\m{GK}$ can be axiomatized either by adding the Fischer Servi axioms for intuitionistic modal logic $\mathsf{IK}$ (see~\cite{FS84b}) to the union of the axioms for both fragments, or by adding the prelinearity axiom for G{\"o}del logic to $\mathsf{IK}$. The fact that the completeness proof strongly relies on assigning intermediate values to the accessibility relation left still open the question of finitely axiomatizing the logic with two modalities of the $\textbf{GK}^c$-models. Indeed, the question of finding an axiomatization complete with respect to the underlying full logic  (with both modal operators) of $\textbf{GK}^c$-models  remains unsolved.

The finite model property with respect to Kripke semantics fails for the bi-modal logics arising from $\textbf{GK}$ and  from $\textbf{GK}^c$. Nevertheless, an alternative semantics for both these logics is introduced in in
\cite{CaMe13}. This is proven to be complete for those logics and it enjoys the finite model property. Further, the size of the model is bounded in terms of the length of the formula under study, which allows the authors to prove decidability of validity in both logics. 

The main contribution of this paper is to establish an axiomatization for the  bi-modal logic arising from the class of $\textbf{GK}^c$-models, closing the open problem and obtaining a full characterization of the main minimal G\"odel modal logics.

The paper is organized as follows: in Section \ref{sec:preliminaries} we introduce some necessary definitions and review some important results that will be used in the paper, both from propositional G\"odel logic and from some known modal extensions. In Section \ref{sec:axiomatization}, we present the axiomatic system $\mathcal{GK}^c$, and prove some technical results about it. In Section \ref{sec:completeness} we show strong completeness of the previous axiomatic system with respect to the bi-modal local logic of the $\textbf{GK}^c$-models, and provide an axiomatization for the global deduction too. We later see in Section \ref{sec:extensions} how to axiomatize the logics arising from the most usual frame restrictions (reflexive, transitive, symmetric, serial and euclidean frames). Lastly, in Section \ref{sec:nonInterdef} we prove that the modal operators are not interdefinable (in any way) in this logic, proving that $\mathcal{GK}^c$ does not coincide with any of the logics of the mono-modal fragments.

\section{Preliminaries}\label{sec:preliminaries}

\subsection{Propositional G\"odel logic}
Let $\mathcal{L}(V)$ be the set of formulas built over a countable set of propositional variables $V$ with the binary symbols $\vee ,\wedge ,\rightarrow$ and constant $\bot$, and where other propositional connectives are defined as usual: $\top \coloneqq \bot \rightarrow \bot, \neg \varphi \coloneqq \varphi \rightarrow \bot, \varphi \leftrightarrow \psi \coloneqq (\varphi \rightarrow \psi) \wedge (\psi \rightarrow \varphi)$.

Let us denote by $\mathcal{G}$ the extension of Hilbert propositional calculus with the prelinearity axiom $(\varphi \rightarrow \psi) \vee (\psi \rightarrow \varphi)$\cite{Dm59}. This system is known to be equivalent, for instance, to Haj\`ek's Basic Logic BL extended with idempotency of the monoidal operation \cite{Ha98}. 
For the sake of self-containment, let us introduce an axiomatization of $\mathcal{G}$:
\begin{align*}
(A1)  &\ \varphi \rightarrow (\psi \rightarrow \varphi) & 
(A8) &\  (\varphi \rightarrow \psi) \rightarrow ((\chi \rightarrow \varphi) \rightarrow (\chi \rightarrow \psi))\\
(A2) &\ (\varphi \wedge \psi) \rightarrow \varphi & 
(A9) &\ (\varphi \rightarrow (\psi \rightarrow \chi)) \rightarrow (\psi \rightarrow (\varphi \rightarrow \chi)) \\
(A3) &\ (\varphi \wedge \psi) \rightarrow \psi & 
(A10) &\  ((\varphi \rightarrow \chi) \wedge (\psi \rightarrow \chi)) \rightarrow ((\varphi \vee \psi) \rightarrow \chi)\\
(A4) &\  \varphi \rightarrow (\psi \rightarrow (\varphi \wedge \psi)) & 
(A11) &\  (\varphi \rightarrow (\psi \rightarrow \chi)) \rightarrow ((\varphi \wedge \psi) \rightarrow \chi) \\
(A5) &\  (\bot \rightarrow \varphi) \wedge (\varphi \rightarrow \top) &
(A12) &\  ((\chi \rightarrow \varphi) \wedge (\chi \rightarrow \psi)) \rightarrow (\chi \rightarrow (\varphi \wedge \psi))\\
(A6) &\   \varphi \rightarrow (\varphi \vee \psi) & 
(A13) &\  (\varphi \rightarrow (\varphi \rightarrow \psi)) \rightarrow (\varphi \rightarrow \psi) \\
(A7) &\  \psi \rightarrow (\varphi \vee \psi)& 
(A14) &\  (\varphi \rightarrow \psi) \vee (\psi \rightarrow \varphi)\\
(MP) &\  \varphi, \varphi \rightarrow \psi \vdash \psi
\end{align*}

$\vdash_{\mathcal{G}}$ denotes the usual deduction in $\mathcal{G}$, and we will write $\vdash_{\mathcal{G}} \varphi$ whenever $\emptyset \vdash_{\mathcal{G}} \varphi$. We will use this notation convention for all other axiomatic systems in the paper.

The equivalent algebraic semantics of $\vdash_{\mathcal{G}}$ is that of the so-called G\"odel algebras, namely semilinear Heyting algebras. This is the variety generated by the \termDef{Standard G\"odel algebra}, the structure
$[0,1]_G = \langle [0,1], \wedge, \vee, \rightarrow, 0\rangle$, where $\wedge$ and $\vee$ are the usual minimum and maximum in $[0,1]$, and  \[ a \rightarrow b = \begin{cases} 1 &\hbox{ if }  a \leq b, \\ b&\hbox{ otherwise}.\end{cases}\]

Since no confusion might arise, we will write, as usual, the same symbols to denote both the syntactic operator in the language and the corresponding operation in the standard G\"odel algebra, for what concerns the previous propositional connectives.
In order to lighten the notation, for any G\"odel homomorphism $h$ and a (possibly infinite) set of formulas $\Gamma$, we shall write $h(\Gamma)$ to denote the set $\{h(\gamma) \colon \gamma \in \Gamma\}$.
Moreover, as usual, for a non-empty finite set of formulas $\Gamma = \{\gamma_1, \ldots, \gamma_n\}$, we write $\bigwedge \Gamma$ to denote the formula $\gamma_1 \wedge \ldots \wedge \gamma_n$, and the analogous for $\bigvee \Gamma$. Further, we use the conventions
\[\bigwedge_{\gamma \in \emptyset} \gamma \coloneqq \top \qquad \bigvee_{\gamma \in \emptyset} \gamma \coloneqq \bot \]
and the analogous for the infimum/supremum of the empty set over elements in $[0,1]$.

Let $\models_{[0,1]_G} $ denote the usual consequence over the standard G\"odel algebra, i.e., for arbitrary $\Gamma \cup \{\varphi\} \subseteq \mathcal{L}(V)$,
\[\Gamma \models_{[0,1]_G} \varphi \text{ iff for all }h \in Hom(\mathcal{L}(V),[0,1]_G), \ h(\Gamma) \subseteq \{1\} \text{ implies } h(\varphi) = 1.\]
 
It is known that $\vdash_{\mathcal{G}}$ does not only enjoy strong completeness with respect to $\models_{[0,1]_G}$, but that this completeness extends to deductions from arbitrary (possibly infinite) theories. Moreover, this implies also an order-preserving completeness that will be useful in the next sections. Let us summarize this facts.
\begin{proposition}\label{prop:propStrongCompleteness}
	Let $\Gamma \cup \{\varphi\} \subseteq \mathcal{L}(V)$\footnote{The premise that $V$ is countable cannot be ignored here.}. The following are equivalent:
	\begin{enumerate}
		\item $\Gamma \vdash_{\mathcal{G}} \varphi$,
		\item $\Gamma \models_{[0,1]_G} \varphi$,
		\item For any $h\in Hom(\mathcal{L}(V),[0,1]_G)$ it holds that $\bigwedge_{\gamma \in \Gamma} h(\gamma) \leq h(\varphi)$.
	\end{enumerate}	
\end{proposition}
\begin{proof}
The equivalence between \textit{1.} and \textit{2.} is proven in \cite{Dm59}, and also see for instance \cite[Th. 4.2.17]{Ha98}. On the other hand, 
	\textit{3.} trivially implies \textit{2.} It can be seen that \textit{2.} implies \textit{3.} easily since any order preserving mapping in $[0,1]$ is a G\"odel endomorphism.

	%
\end{proof}

Moreover, $\vdash_{\mathcal{G}}$ enjoys the usual Deduction Theorem (D.T), i.e., for any $\Gamma \cup \{\psi, \varphi\} \subseteq \mathcal{L}$, \[\Gamma, \psi \vdash_{\mathcal{G}} \varphi \text{ if and only if } \Gamma \vdash_{\mathcal{G}} \psi \rightarrow \varphi.\]

Before continuing, let us exhibit a formula that is valid in $\mathcal{G}$ which will be used in Section \ref{sec:completeness}.

\begin{lemma}\label{lemma:theoremsplitting}
	$\vdash_{\mathcal{G}} (((\chi \rightarrow \phi) \rightarrow \phi) \wedge (\phi \rightarrow \psi)) \vee (((\chi \rightarrow \phi) \rightarrow \phi) \rightarrow (\psi \rightarrow \phi))$.
\end{lemma}
\begin{proof}
	Let $h \in Hom(\mathcal{L}, [0,1]_G)$ and suppose $h(((\chi \rightarrow \phi) \rightarrow \phi) \wedge (\phi \rightarrow \psi)) < 1$. Thus, either $h((\chi \rightarrow \phi) \rightarrow \phi) < 1$ or $h(\phi \rightarrow \psi) < 1$. In the second case, by prelinearity, $h(\psi \rightarrow \phi) = 1$, and so, $h(((\chi \rightarrow \phi) \rightarrow \phi) \rightarrow (\psi \rightarrow \phi)) = 1$. 
	On the other hand, if $h((\chi \rightarrow \phi) \rightarrow \phi) < 1$, due to the definition of the implication in $[0,1]_G$,  necessarily $h((\chi \rightarrow \phi) \rightarrow \phi) = h(\phi)$. Then,
	$h(((\chi \rightarrow \phi) \rightarrow \phi) \rightarrow (\psi \rightarrow \phi)) = h(\phi) \rightarrow h(\psi \rightarrow \phi) = 1$. We conclude the proof relying in completeness of $\vdash_{\mathcal{G}}$ (Proposition \ref{prop:propStrongCompleteness}).
\end{proof}

\subsection{G\"odel modal logics}
Let us consider a modal expansion of G\"odel logic with two operators $\Box$ and $\Diamond$.  The set of formulas $\mathcal{L}_{\Box \Diamond }(V)$ 
is built as $\mathcal{L}(V)$ (always assuming countability of the set of propositional variables $V$) but extending the set operations with two unary symbols $\Box $ and $\Diamond$. Whenever $V$ is clear from the context we will simply write $\mathcal{L}_{\Box \Diamond}$.
We will sometimes refer to the mono-modal expansions of G\"odel logic, namely, those extending the propositional language either with only $\Box$ ($\mathcal{L}_\Box(V)$) or only $\Diamond$ ($\mathcal{L}_\Diamond(V)$).  We will say that $\varphi$ is a mono-modal formula whenever $\varphi \in \mathcal{L}_\Box(V) \cup \mathcal{L}_\Diamond(V)$.

In the style introduced by Fitting \cite{Fi92a,Fi92b} and studied in the works mentioned in the introduction, we define the G\"odel Modal Logic as arising from its semantic definition. This is given by enriching usual Kripke models with evaluations over the previous standard algebra, as in \cite{CaRo10, CaRo15} and others. Formally:

\begin{definition}
	A \termDef{G\"odel-Kripke model} $\m{M}$  is a structure $\langle W, R, e\rangle$ where $W$ is a non-empty set of so-called worlds, and $R \colon W \times W \rightarrow[0,1]$ and $e \colon V \times W \rightarrow [0,1]$ are arbitrary mappings.
	
	Whenever $R \colon W \times W \rightarrow \{0,1\}$ we will say that the model is \termDef{crisp}, and write $Rvw$ to denote $R(v,w) = 1$.
\end{definition}

The evaluation $e$ can be uniquely extended to a map with domain $W \times \mathcal{L}_{\Box\Diamond}$ in such a way that it is a propositional G\"odel homomorphism (for the propositional connectives) and where the modal operators are interpreted as infima and suprema Mostowski style \cite{Mo57}, namely:
\begin{itemize}
	\item $e(v, \bot) \coloneqq 0$,
	\item $e(v, \varphi \star \psi) \coloneqq e(v, \varphi) \star e(v, \psi)$  for $\star \in \{\wedge, \vee, \rightarrow\}$,
	\item $e(v,\Box \varphi) \coloneqq \bigwedge_{w \in W}( R(v,w) \rightarrow e(w, \varphi)),$
	\item $e(v,\Diamond \varphi) \coloneqq \bigvee_{w \in W}( R(v,w) \wedge e(w, \varphi)).$
\end{itemize}

Truth and logical {
entailment over the whole class of models, and over the crisp ones, are defined as follows. Observe that truth and entailment from G\"odel propositional logic is world-wise preserved.
\begin{definition}
\leavevmode
\makeatletter
\@nobreaktrue
\makeatother
\begin{itemize}
		\item Formula $\varphi$ \termDef{is true at world }$v$ in the model $\m{M}$, and write $\m{M}, v \models \varphi$ if and only if $e(v, \varphi) = 1$. Formula $\varphi$ \termDef{is true at model }$\m{M}$, and write $\m{M} \models \varphi$ if and only if $\m{M}, v \models \varphi$ for all $v \in W$.
		\item We say that the formula $\varphi$ \termDef{follows locally} from the set of formulas $\Gamma$, and write  $\Gamma \models_{\m{GK}^c}  \varphi$  $(\Gamma \models_{\m{GK}} \varphi)$ if and only if for any crisp G\"odel-Kripke model (G\"odel-Kripke model) $\m{M}$,
		\[\text{ for all }v \in W:\quad  \m{M}, v \models \gamma\text{ for all } \gamma \in \Gamma \quad \text{ implies } \quad \m{M}, v \models \varphi\]
		\item We say that the formula $\varphi$ \termDef{follows globally} from the set of formulas $\Gamma$, and write $\Gamma \models^g_{\m{GK}^c} \varphi$ $(\Gamma \models^g_{\m{GK}} \varphi)$ if and only if for any crisp G\"odel-Kripke model (G\"odel-Kripke model) $\m{M}$,
		\[\m{M} \models \gamma\text{ for all } \gamma \in \Gamma \quad  \text{ implies} \quad \m{M} \models \varphi\]
	\end{itemize}
\end{definition}

Observe that the set of theorems of the local and the global logics are clearly the same, but the deduction systems, as it happens in the classical case, are not: $\varphi \models^g_{\m{GK}} \Box \varphi$, but that is not the case in the local consequence.
Along this work, we will mainly study the local deduction, and we will prove at the end of Section \ref{sec:completeness} a completeness result for the global logic, building on the completeness of the local one.

In \cite{CaRo15} the authors study the logic $\models_{\m{GK}}$ defined in the above way with both $\Box$ and $\Diamond$, and the axiomatic system $\mathcal{GK}$ is introduced and proven complete.
$\mathcal{GK}$ in \cite{CaRo15} is defined as the extension of the Intuitionistic modal logic IK by Fischer-Servi (see eg. \cite{FS84b}) with the prelinearity axiom, in the same fashion that it is done in the propositional case. This coincides with the system resulting from extending the calculus of G\"odel-Dummet propositional logic $\mathcal{G}$
by the following set of axioms and rules:
\begin{align*}
(K_\Box) &\  \Box(\varphi \to \psi) \to (\Box \varphi \to \Box \psi) & (K_\Diamond)&\  \Diamond(\varphi \lor \psi) \to (\Diamond\varphi  \lor \Diamond \psi)\\
(FS1) &\   \Diamond(\varphi \rightarrow \psi) \rightarrow (\Box \varphi \rightarrow \Diamond \psi) & (FS2) &\  (\Diamond\varphi  \to \Box \psi) \to \Box(\varphi \to \psi)\\
(F_\Diamond) &\  \neg \Diamond \bot & \\
(N_\Box) &\  \vdash \varphi \text{ implies }\vdash \Box \varphi & (N_\Diamond)&\  \vdash \varphi \rightarrow \psi \text{ implies } \vdash \Diamond \varphi \rightarrow \Diamond \psi
\end{align*}

\begin{theorem}[Th. 3.1, \cite{CaRo15}]\label{th:compGK}
Let $\Gamma, \varphi \subseteq \mathcal{L}_{\Box\Diamond}(V)$. Then
\[\Gamma \vdash_{\mathcal{GK}} \varphi \text{ if and only if } \Gamma \models_{\m{GK}} \varphi.\]
\end{theorem}

In  \cite{CaRo15} it is pointed out that an alternative axiomatization of the previous system can be given by replacing $(FS1)$ with the axiom scheme
\begin{align*}
(P) &\ \Box(\varphi \to \psi) \to (\Diamond \varphi \to \Diamond \psi)
\end{align*}
and removing the rule $(N_\Diamond)$. In general, we will be using this second presentation of the logic $\mathcal{GK}$, particularly when facing a proof by induction on the length of a derivation in the logic.

Some formulas valid in $\mathcal{GK}$ that will be used below are the following:
\begin{align*}
(T1) &\  \Box(\varphi \wedge \psi) \leftrightarrow \Box \varphi  \wedge \Box \psi\\
 (T2)&\  ((\Box \varphi \rightarrow \Diamond \psi) \rightarrow \Diamond \psi) \rightarrow \Box((\varphi \rightarrow \psi) \rightarrow \psi) \vee \Diamond \psi
\end{align*}
It is easy to check both of them are valid in $\models_{\m{GK}}$, and so theorems of $\mathcal{GK}$.

In \cite{CaRo15}, the axiomatization of the logic arising from $\models_{\m{GK}^c}$ in the language with two modalities is left as an open problem. On the other hand, the corresponding mono-modal fragments have been studied and axiomatized in \cite{CaRo10} and \cite{MeOl11}. 

The axiomatic system $\mathcal{GK}^c_\Box$ is introduced in \cite{CaRo10} (under the name of $\mathcal{G}_\Box$). It is the extension of $\mathcal{G}$ with the following axiom schemata and rule:
\begin{align*}
(K_\Box) &\  \Box(\varphi \to \psi) \to (\Box \varphi \to \Box \psi) & (N_\Box) &\   \vdash \varphi \text{ implies }\vdash \Box \varphi \\
(Z_\Box) &\  \neg \neg \Box \varphi \rightarrow \Box \neg \neg \varphi
\end{align*}

\begin{theorem}[Th. 4.2, \cite{CaRo10}]\label{th:compGKBox}
Let $\Gamma, \varphi \subseteq \mathcal{L}_{\Box}(V)$. Then 
\[\Gamma \vdash_{\mathcal{GK}^c_\Box} \varphi \text{ if and only if } \Gamma \models_{\m{GK}^c} \varphi.\]
\end{theorem}

On the other hand, the $\Diamond$-fragment is studied in \cite{MeOl11}. 
The axiomatic system $\mathcal{GK}^c_\Diamond$ is introduced there (under the name of $HGK_\Diamond$). It is the extension of $\mathcal{G}$ with the following axiom schemata and rule:
\begin{align*}
(Z_\Diamond) &\  \Diamond \neg \neg \varphi \rightarrow \neg \neg \Diamond \varphi & (K_\Diamond) &\   \Diamond(\varphi \vee \psi) \rightarrow (\Diamond \varphi \vee \Diamond \psi) \\
(F_\Diamond) &\  \neg \Diamond \bot &
(Nec_\Diamond) &\  \vdash (\varphi \rightarrow \psi) \vee \chi \text{ infer } \vdash(\Diamond \varphi \rightarrow \Diamond \psi) \vee \Diamond \chi
\end{align*}

\begin{theorem}[Th. 5.8, \cite{MeOl11}] \label{th:compGKDiamond}
Let $\Gamma, \varphi \subseteq \mathcal{L}_{\Diamond}(V)$. Then
\[\Gamma \vdash_{\mathcal{GK}_\Diamond^c} \varphi \text{ if and only if } \Gamma \models_{\m{GK}^c} \varphi.\]
\end{theorem}

\section{The logic $\mathcal{GK}^c$}\label{sec:axiomatization}

As we said before, the axiomatization of the logic $\models_{\m{GK}^c}$ with both $\Box$ and $\Diamond$ is an open problem. It is also not known whether the system $\mathcal{GK} \cup \mathcal{GK}^c_\Box \cup \mathcal{GK}^c_\Diamond$ might be complete with respect to $\models_{\m{GK}^c}$ for sets of formulas in $\mathcal{L}_{\Box\Diamond}$.

We propose here the axiomatic system $\mathcal{GK}^c$, which extends $\mathcal{GK}$ with an axiom that is also used in the field of positive modal logics (see for instance \cite{Du95}). We will prove that this system is strongly complete with respect to $\models_{\m{GK}^c}$, solving the open problem stated above.

\begin{definition}

The logic $\mathcal{GK}^c$ is defined by adding to  $\mathcal{GK}$ the following axiom scheme
\begin{align*}
(Cr) &\  \Box(\varphi \vee \psi) \rightarrow (\Box \varphi \vee \Diamond \psi)
\end{align*}
\end{definition}

Let us pay some attention to the relation of $\mathcal{GK}^c$ with respect to the existing axiomatizations of the mono-modal fragments. 
Since in \cite{CaRo15} it is proven that $\vdash_{\mathcal{GK}} (Z_\Box)$, and the other axiom schemata and rule from $\mathcal{GK}_\Box^c$ are explicitly included in the definition of $\mathcal{GK}$, we get the following.
\begin{remark}\label{obs:weakerBox}
For $\Gamma, \varphi \subseteq \mathcal{L}_\Box(V)$, 
$\Gamma \vdash_{\mathcal{GK}^c_\Box} \varphi$ implies that $\Gamma \vdash_{\mathcal{GK}^c} \varphi$.
\end{remark}

 The same relation of $\mathcal{GK}^c_\Diamond$ with respect to $\mathcal{GK}^c$ can be proven too.
\begin{lemma}\label{lem:weakerDiamond}
For $\Gamma, \varphi \subseteq \mathcal{L}_\Diamond(V)$, 
$\Gamma \vdash_{\mathcal{GK}^c_\Diamond} \varphi$ implies that $\Gamma \vdash_{\mathcal{GK}^c} \varphi$.
\end{lemma}
\begin{proof}
It is only needed to prove that the formula $(Z_\Diamond)$ and the rule $(Nec_\Diamond)$ can be derived in $\mathcal{GK}^c$. 

Concerning $(Z_\Diamond)\ \Diamond \neg \neg \varphi \rightarrow \neg \neg \Diamond \varphi$, observe that $\vdash_{\mathcal{GK}^c}\Diamond \neg \neg \varphi \rightarrow (\Box \neg \varphi \rightarrow \Diamond \bot)$, by $(FS1)$, and so, from $(F_\Diamond)$ it follows that 
\begin{equation}
\vdash_{\mathcal{GK}^c} \Diamond \neg \neg \varphi \rightarrow \neg \Box \neg \varphi
\end{equation} 

Further, $\vdash_{\mathcal{GK}^c} \neg \Diamond \varphi \rightarrow ((\Diamond \varphi) \rightarrow \Box \bot) $ (since $\vdash_{\mathcal{GK}^c} \bot \rightarrow \Box \bot)$, and thus, applying $(FS2)$, we know that
$\vdash_{\mathcal{GK}^c} \neg \Diamond \varphi \rightarrow \Box \neg \varphi$.
Thus $\vdash_{\mathcal{GK}^c} \neg \Box \neg \varphi \rightarrow \neg \neg \Diamond \varphi$, and using the implication (1) proven above we conclude
$\vdash_{\mathcal{GK}^c} \Diamond \neg \neg \varphi \rightarrow \neg \neg \Diamond \varphi$. 

Let us now prove that 
the inference rule
\[(R_\Diamond) \quad \vdash \varphi \vee (\psi \rightarrow \chi) \text{ implies } \vdash \Diamond \varphi \vee (\Diamond \psi \rightarrow \Diamond \chi),\]
is also derivable in $\mathcal{GK}^c$. 
This can be proven by first applying rule $(N_\Box)$ to the premise, getting $\Box (\varphi \vee (\psi \rightarrow \chi) )$. By $(Cr)$ and M.P. it follows that $\Diamond \varphi \vee \Box(\psi \rightarrow \chi)$. Applying $(P)$ to the second part of the previous disjunction, we reach the conclusion.
\end{proof}

\begin{corollary}
Let $\varphi$ be a mono-modal formula. Then 
$\models_{\m{GK}^c} \varphi$ implies $\vdash_{{\mathcal{GK}^c}} \varphi$. 
\end{corollary}
\begin{proof}
If $\varphi$ is a formula using only the $\mathtt{M}$ modality (for $\mathtt{M}\in \{\Box, \Diamond\}$), we know that $\models_{\m{GK}^c} \varphi$ implies that $\vdash_{\mathcal{GK}^c_{\mathtt{M}}} \varphi$ (Theorems \ref{th:compGKBox} and \ref{th:compGKDiamond}).  From Remark \ref{obs:weakerBox} and Lemma \ref{lem:weakerDiamond}, then also 
$\vdash_{\mathcal{GK}^c} \varphi$.
\end{proof}

Using this, let us exhibit some additional valid formulas from $\mathcal{GK}^c$ that will be useful in the next section.

	\begin{lemma}
The following formulas are provable in $\mathcal{GK}^c$:	
	\begin{align*}
	(T3) &\  (\Box\varphi \to \Diamond \varphi) \vee \Box \bot,\\
	(T^<_{\Box}) &\  ((\Box \psi \rightarrow \Box \varphi) \rightarrow \Box \varphi) \rightarrow ((\Box((\psi \rightarrow \varphi)\rightarrow \varphi) \rightarrow \Box \varphi) \rightarrow \Box \varphi),\\
(T^<_{\Diamond}) &\  ((\Diamond \psi \rightarrow \Diamond \varphi) \rightarrow \Diamond \varphi)  \rightarrow \Diamond((\psi \rightarrow \varphi) \rightarrow \varphi)
\end{align*}
\end{lemma}

	\begin{proof}
	$(T^<_{\Box})$ and $(T^<_{\Diamond})$ are mono-modal formulas, and they are easy to check in $\models_{\m{GK}^c}$. Then, from the previous corollary, we get they are derivable in $\mathcal{GK}^c$ too.

	$(T3)$ follows easily from $(Cr)$ and the fact that $\mathcal{GK}^c$ extends $\mathcal{GK}^c_\Box$. 
	Indeed, since $\varphi \rightarrow \bot \vee \varphi $ is a theorem of $\mathcal{G}$, applying $N_\Box$ and subsequently $K$ axiom to it, we get that $\vdash_{\mathcal{GK}^c} \Box \varphi \rightarrow \Box (\bot \vee \varphi)$. Now, by axiom $(Cr)$ and transitivity of the implication, it follows that 
	$\vdash_{\mathcal{GK}^c} \Box \varphi \rightarrow (\Box \bot  \vee \Diamond \varphi)$, and by distributivity of $\rightarrow$ over $\vee$, $\vdash_{\mathcal{GK}^c} (\Box \varphi \rightarrow \Box \bot) \vee  (\Box \varphi  \rightarrow \Diamond \varphi)$.
	
	From here, using that $\vdash_{\mathcal{GK}^c_\Box} (\Box \varphi \rightarrow \Box \bot) \rightarrow \neg \Box \varphi \vee \Box \bot$, and that 
	$\vdash_{\mathcal{G}} \neg \chi_1 \rightarrow (\chi_1 \rightarrow \chi_2)$ for any $\chi_1, \chi_2$ (and so, 
	$\vdash_{\mathcal{GK}^c} \neg \Box \varphi \rightarrow (\Box \varphi \rightarrow \Diamond \varphi)$), we conclude $(T3)$. 
		\end{proof}

We denote by  $\texttt{M}V\coloneqq \{\Box \theta ,\Diamond \theta \colon \theta \in \mathcal{L}_{\Box
\Diamond }(V)\}$, the set of formulas in $\mathcal{L}_{\Box \Diamond }(V)$ starting with a modal symbol, $\Box$ or $\Diamond$. 
If we use this set as names for fresh variables (i.e., not in $V$), clearly 
$\mathcal{L}_{\Box \Diamond }(V) = \mathcal{L}(V \cup \texttt{M}V)$ as sets\footnote{They are not the same if seen as the respective formula algebras, since they have different types - the first one has more operations.}. 
That is to say, any formula in $\mathcal{L}_{\Box\Diamond }(V)$ may be seen as a propositional G\"odel  formula built from the extended set
of  propositional variables $V\cup \texttt{M}V$. 
This allows us to abuse the definition of Homomorphism (which technically is given only for two algebras of the same type), and write 
$Hom(\mathcal{L}_{\Box \Diamond }(V), [0,1]_G)$ to denote $Hom(\mathcal{L}(V \cup \texttt{M}V), [0,1]_G)$. This syntactic 
 association allows us to take advantage of Lemma \ref{lemma:propReduction}.

Let us denote by $Th(\mathcal{GK}^c)$ the set of theorems of $\mathcal{GK}^c$, i.e., the formulas that can be derived in $\mathcal{GK}^c$ from the empty set. 

It is easy to see that deductions in $\mathcal{GK}^c$ can be reduced to derivations in pure propositional G\"odel logic $\vdash_{\mathcal{G}}$ with a certain set of premises. The proof follows immediately from the 
fact that the only non-propositional inference rule from $\mathcal{GK}^c$ is restricted to the set of theorems. 

\begin{lemma}\label{lemma:propReduction}
For any $\Gamma \cup \{\varphi\} \subseteq \mathcal{L}_{\Box\Diamond}$,
 \[\Gamma \vdash_{\mathcal{GK}^c} \varphi \text{ if and only if } Th(\mathcal{GK}^c),  \Gamma \vdash_{\mathcal{G}} \varphi.\]
\end{lemma}

\noindent It is also easy to see that $\vdash_{\mathcal{GK}^c}$ still enjoys the D.T, namely for any $\Gamma \cup \{\psi, \varphi\} \subseteq \mathcal{L}_{\Box \Diamond}$,
\begin{equation}
\tag{D.T}
\Gamma, \psi  \vdash_{\mathcal{GK}^c} \varphi \text{\textit{ if and only if }} \Gamma  \vdash_{\mathcal{GK}^c} \psi \rightarrow \varphi
\end{equation}

In addition, we can prove the following meta-rule, which will be useful to prove completeness of $\vdash_{\mathcal{GK}^c}$. As usual, for an arbitrary set $\Gamma \subseteq \mathcal{L}_{\Box\Diamond}$, we let \[\Box \Gamma \coloneqq \{\Box \gamma\colon \gamma \in \Gamma\}.\]

\begin{lemma}\label{lemma:metarules}
	For any $\Gamma \cup \{\varphi\} \subseteq \mathcal{L}_{\Box \Diamond}$,
	\begin{equation}
	\tag{$M_\Box$}
	\Gamma \vdash_{\mathcal{GK}^c} \varphi \text{ implies } \Box \Gamma \vdash_{\mathcal{GK}^c} \Box \varphi
	\end{equation}
\end{lemma}

\begin{proof}
We reason by  induction on the length of the derivation of $\varphi$ from $\Gamma$ in $\mathcal{GK}^c$. We use the presentation of $\mathcal{GK}^c$ with only inference rules M.P (from $\mathcal{G}$) and $(N_{\Box})$.

 If $\vdash_{\mathcal{GK}^c} \varphi$, then the step follows by the necessitation rule. 
 Otherwise $\Gamma \vdash_{\mathcal{GK}^c} \chi$ and $\Gamma \vdash_{\mathcal{GK}^c} \chi \rightarrow \varphi$ (since M.P. is the only inference rule affecting not only theorems of the logic). By I.H. $\Box \Gamma \vdash_{\mathcal{GK}^c} \Box \chi$ and $\Box \Gamma \vdash_{\mathcal{GK}^c} \Box(\chi \rightarrow \varphi)$. Applying $K_\Box$ axiom and later MP we get $\Box \Gamma \vdash_{\mathcal{GK}^c} \Box \varphi$.
 	\end{proof}

\section{Completeness of $\mathcal{GK}^c$}\label{sec:completeness}

In this section we will show that $\mathcal{GK}^c$  is complete with respect to the local deduction in $\m{GK}^c$. We will begin by detailing the proof for valid formulas, and at the end of the section we will see that this easily extends to all deductions in the logic. We will also see how we can use this completeness to provide an axiomatization of the global deduction over the same class of models.

For any formula $\rho$  we denote by $\langle \rho \rangle \subseteq \mathcal{L}_{\Box \Diamond }$ the set of subformulas of $\rho$ containing in addition, constants $\bot $ and $\top$.

For each formula $\rho \in \mathcal{L}_{\Box \Diamond}$ that is not a theorem of $\mathcal{GK}^c$, we will build a crisp G\"odel-Kripke model $\m{M}^{\rho}$ where there is indeed a world in which $\rho$ is evaluated to less than $1$.
In order to do so, we will define a structure in a similar fashion to the canonical model from \cite{CaRo15}, and we will see it is canonical for $\langle \rho \rangle$.\footnote{Meaning that each world of the model is a G\"odel homomorphism satisfying all theorems of the modal logic, and at each world $h$, and for each formula $\psi \in \langle \rho \rangle$, $e(h, \psi) = h(\psi)$, taking into account the syntactic convention of the set equality $\mathcal{L}_{\Box\Diamond}(V) = \mathcal{L}(V \cup \mathtt{M}V)$.}

The \emph{canonical model} $\m{M}^{\rho}=\langle W^{\rho},R^{\rho},e^{\rho}\rangle$ is defined as follows:
\begin{itemize}
\item $W^{\rho}$ is the set $\{u \in Hom(\mathcal{L}_{\Box\Diamond}(V), [0,1]_G)\colon  u(Th(\mathcal{GK}^c))\subseteq\{1\} \}$.\footnote{Recall this notation stands for $Hom(\mathcal{L}(V \cup \texttt{M}V), [0,1]_G)$.}
\item $R^{\rho}w u$ if and only if $\forall \psi \in \langle \rho \rangle : w(\Box \psi) \leq u(\psi) \mbox{ and } u(\psi) \leq w(\Diamond\psi).$

\item $e^{\rho}(u,p)=u(p)$ for any $p\in V$.
\end{itemize}

The previous structure is, by definition, a G\"odel-Kripke model. The main idea behind the definition is that, if $\rho \not \in Th(\mathcal{GK}^c)$, then from Lemma \ref{lemma:propReduction}  and strong standard completeness of $\vdash_{\mathcal{G}}$ (Proposition \ref{prop:propStrongCompleteness}), there is $h \in Hom(\mathcal{L}({V\cup \texttt{M}V}), [0,1]_G)$ such that $h(Th(\mathcal{GK}^c)) \subseteq \{1\}$ and $h(\rho) < 1$.
To use this homomorphism in proving that this model is indeed a counter-model for $\rho$, we need to see that $e(h, \rho) = h(\rho)$.
We will do so by proving a version of the usual Truth-Lemma relative to $\langle \rho \rangle$, which can be done because of the way we defined $R^{\rho}$ above.

Let us introduce some notation to simplify  the reading of the results below. For $u \in W^{\rho}$, $\alpha \in [0,1]$, modality $\texttt{M} \in \{\Box, \Diamond\}$ and $\triangledown \in \{<,>,=\}$  put
\[\texttt{M}^{\triangledown \alpha}_u \coloneqq \{\psi \in \langle \rho \rangle\colon u(\texttt{M}\psi) \triangledown \alpha\}\]

Moreover, we will denote the versions of the above sets not restricted to formulas in $\langle \rho \rangle$ by ${^\ast\texttt{M}}^{\triangledown \alpha}_u$. In fact, we will be only using one of these sets, namely
\[{^\ast\Box}^{= 1}_u = \{\psi \in \mathcal{L}_{\Box\Diamond}(V) \colon u(\Box \psi) = 1\}.\]
A trivial observation about the above sets is that 
for any $\psi  \in \texttt{M}^{\triangledown \alpha}_u$, $u(\texttt{M}\psi) \triangledown \alpha$.

We will sometimes refer to the formulas $\bigwedge U$ or $\bigvee U$ for some of the above $\texttt{M}^{\triangledown \alpha}_u $ sets 
(since they are always finite, this is well defined, see the preliminaries section). Recall that, by convention, we assume that if $U=\emptyset$, these are respectively the formulas $\top$ and $\bot$. 

Let us begin by proving some results that will later allow to give an easy proof of the Truth-Lemma for the $\Box$-formulas.

\begin{lemma}\label{lemma:deltaalpha}
	Let $\alpha < 1$ and $\varphi  \in\Box^{=\alpha}_u$. Let \footnote{While $\delta$ depends on $u$ and $\varphi$, we have chosen to omit these elements from the name of the formula, since they are clear from the context and the notation gets much heavier if we use $\delta_u^\varphi$.}
	\[\delta \coloneqq (\bigwedge \Box^{>\alpha}_u \rightarrow \varphi) \rightarrow \varphi\]
	Then $u(\Box \delta) > \alpha$.
	
\end{lemma}
\begin{proof}

	From $(T1)$ -namely, distributivity of $\wedge$ and $\Box$-  we know that
	$u(\Box \bigwedge \Box^{>\alpha}_u) = u(\bigwedge \Box \Box^{>\alpha}_u)$. Since for any $\psi \in \Box^{>\alpha}_u$ 
	by definition $u(\Box \psi) > \alpha$,  we get that $u(\Box(\bigwedge \Box^{>\alpha}_u))  > \alpha$, and in particular, since $\varphi \in \Box^{=\alpha}_u$, 
	$u(\Box \bigwedge \Box^{>\alpha}_u)  > u(\Box \varphi)$. Then, from the characteristics of G\"odel implication, it follows that
	\begin{equation*}
	u((\Box \bigwedge \Box^{>\alpha}_u \rightarrow \Box \varphi) \rightarrow \Box \varphi) = 1.
	\end{equation*}
	
Consider now the formula $(T^<_\Box)$, valid in $\mathcal{GK}^c$. We can substitute in its premise the previous formula, and by M.P. we know that 
	\begin{equation*}
u((\Box (( \bigwedge \Box^{>\alpha}_u \rightarrow \varphi)\rightarrow \varphi) \rightarrow \Box \varphi) \rightarrow \Box \varphi) = 1.
	\end{equation*}
	
	From the definition of G\"odel implication, and since $u(\Box \varphi) < 1$, the above implies that
 $u(\Box \delta) = u(\Box((\bigwedge \Box^{>\alpha}_u  \rightarrow \varphi)\rightarrow \varphi)) > u(\Box \varphi) = \alpha$, concluding the proof.
\end{proof}

The next remark is a matter of expanding the definitions.
\begin{remark}\label{obs:deltacond}
For any G\"odel homomorphism $v$, if $v(\delta) = 1$ and $v(\varphi) < 1$ then $v(\varphi) < v(\psi)$ for all $\psi \in \Box^{>\alpha}_u$.
\end{remark}

\begin{proposition} \label{prop:witnessBox}
	Let $\alpha < 1$ and $\varphi  \in\Box^{=\alpha}_u$. Then there exists $h \in Hom(\mathcal{L}_{\Box \Diamond}(V), [0,1]_G)$ such that
	\begin{enumerate}[align=left]
		\item[(C1)] $h(Th(\mathcal{GK}^c)) \subseteq\{ 1\}$,
		\item[(C2)] $h({^\ast\Box}^{=1}_u) \subseteq \{1\}$,
		\item[(C3)] $h(\psi) < 1$ for all $\psi \in \Diamond^{<1}_u$,
		\item[(C4)] $h(\varphi) < h(\psi)$ for all $\psi \in \Box^{>\alpha}_u$.
	\end{enumerate}
\end{proposition}
\begin{proof}
Recall from Lemma \ref{lemma:theoremsplitting} that for any $\chi, \phi, \psi$, it holds that 
	\[\vdash_{\mathcal{G}} (((\chi \rightarrow \phi) \rightarrow \phi) \wedge (\phi \rightarrow \psi)) \vee (((\chi \rightarrow \phi) \rightarrow \phi) \rightarrow (\psi \rightarrow \phi)).\] 
Substituting  $\chi$ by $\bigwedge \Box^{>\alpha}_u$, $\phi$ by $\varphi$ and $\psi$ by $\bigvee \Diamond^{<1}_u$, and using the wrapping $\delta$ introduced in Lemma \ref{lemma:deltaalpha} we get that 
	\[\vdash_{\mathcal{GK}^c} (\delta \wedge (\varphi \rightarrow \bigvee \Diamond^{<1}_u)) \vee
	 (\delta \rightarrow (\bigvee \Diamond^{<1}_u \rightarrow \varphi)).\]
	
	Applying commutativity of $\vee$, the $(N_\Box)$ rule and axiom $(Cr)$, we get
		\[\vdash_{\mathcal{GK}^c}\Diamond (\delta \wedge (\varphi \rightarrow \bigvee \Diamond^{<1}_u)) \vee \Box(\delta \rightarrow (\bigvee \Diamond^{<1}_u \rightarrow \varphi)).\]
Since $u \in W^{\rho}$, it evaluates the previous formula to $1$ necessarily, and so there are two possible cases:
\begin{itemize}[align=left]
	\item[(A)] Either $u(\Diamond (\delta \wedge (\varphi \rightarrow \bigvee \Diamond^{<1}_u))) = 1$, or
	\item[(B)] $u(\Box(\delta \rightarrow (\bigvee \Diamond^{<1}_u \rightarrow \varphi))) = 1$.
\end{itemize}
We will show that in either case the Proposition can be proven.

\noindent (A) Assume $u(\Diamond (\delta \wedge (\varphi \rightarrow \bigvee \Diamond^{<1}_u))) = 1$, and 
			let us prove
			\begin{equation}\label{eq:B1}
			Th(\mathcal{GK}^c), {^\ast\Box}^{=1}_u, \delta \not \models_{[0,1]_G} (\varphi \rightarrow \bigvee \Diamond^{<1}_u) \rightarrow \bigvee \Diamond^{<1}_u
			\end{equation}
			Suppose the contrary, with a view to contradiction.		
			Using (strong) completeness of $\models_{[0,1]_G}$ with respect to $\mathcal{G}$ (Proposition \ref{prop:propStrongCompleteness}), and then Lemma \ref{lemma:propReduction} (which allows us to move between propositional and modal deductions) and the D.T., it follows that
			\[{^\ast\Box}^{=1}_u \vdash_{\mathcal{GK}^c} (\delta \wedge (\varphi \rightarrow \bigvee \Diamond^{<1}_u)) \rightarrow \bigvee \Diamond^{<1}_u.\]
Applying the meta-rule $M_\Box$ (Lemma \ref{lemma:metarules}), and axioms  $(P)$ and $(K_\Diamond)$ it follows that
			\[\Box{^\ast\Box}^{=1}_u \vdash_{\mathcal{GK}^c} \Diamond (\delta \wedge (\varphi \rightarrow \bigvee \Diamond^{<1}_u)) \rightarrow \bigvee \Diamond \Diamond^{<1}_u.\]
			Going back to propositional (via Lemma \ref{lemma:propReduction} again), it follows that
			\[Th(\mathcal{GK}^c), \Box{^\ast\Box}^{=1}_u \models_{[0,1]_G} \Diamond (\delta \wedge (\varphi \rightarrow \bigvee \Diamond^{<1}_u)) \rightarrow\bigvee \Diamond \Diamond^{<1}_u.\]
			However, this leads to a contradiction, since we can prove $u$ refutes this derivation:
			\begin{itemize}
				\item $u(Th(\mathcal{GK}^c)) \subseteq \{1\}$ (since $u \in W$), and $u(\Box {^\ast\Box}^{=1}_u) \subseteq \{1\}$ by definition. Thus,  the premises of the derivation are met by homomorphism $u$. However, 
				\item $u(\Diamond (\delta \wedge (\varphi \rightarrow \bigvee \Diamond^{<1}_u))) = 1$, 
				since we assumed (A) at the beginning of this part of the proof,  and $u(\bigvee \Diamond\Diamond^{<1}_u) < 1$ by definition. 
				Thus, $u(\Diamond (\delta \wedge (\varphi \rightarrow \bigvee \Diamond^{<1}_u)) \rightarrow \bigvee \Diamond\Diamond^{<1}_u) < 1$, 
				meaning that the conclusion is not satisfied by $u$ and so contradicting the definition of $\models_ {[0,1]_G}$.
			\end{itemize}
			
			This concludes the proof of condition (\ref{eq:B1}). Thus, 
			there exists an homomorphism $h \in Hom(\mathcal{L}_{\Box\Diamond}(V), [0,1]_G)$ that sends the premises of (\ref{eq:B1}) to $1$ and the conclusion to some value strictly less than $1$. We claim this homomorphism $h$ meets the four conditions stated in the Proposition, since: 
			\begin{itemize}
				\item The premises on (\ref{eq:B1}) are sent to $1$ by $h$, so
				$h(Th(\mathcal{GK}^c))  \subseteq \{1\}$ proving \textit{(C1)}, and $h({^\ast\Box}^{=1}_u)  \subseteq \{1\}$ proving  \textit{(C2)}. 
				\item $h((\varphi \rightarrow \bigvee \Diamond^{<1}_u) \rightarrow \bigvee \Diamond^{<1}_u) < 1$ implies that 
				$h(\varphi \rightarrow \bigvee \Diamond^{<1}_u) > h(\bigvee \Diamond^{<1}_u)$. Thus, necessarily, $h(\bigvee \Diamond^{<1}_u) < 1$, proving $h$ satisfies \textit{(C3)}. Further, by the definition of G\"odel implication, it also follows that 
				$h(\varphi) \leq h(\bigvee \Diamond^{<1}_u)$, proving that also $h(\varphi) < 1$.
				
				\item Using again that the premises of (\ref{eq:B1}) are sent to $1$ by $h$, we know that $h(\delta) = 1$. Together with $h(\varphi) < 1$ (from the previous point) and Remark \ref{obs:deltacond}, these imply that $h(\varphi) < h(\psi)$ for any $\psi \in \Box_u^{>\alpha}$, namely, \textit{(C4)}.
			\end{itemize}

\noindent (B) Assume (A) does not hold, and so, (B) is the case, i.e.,
 $u(\Box(\delta \rightarrow (\bigvee \Diamond^{<1}_u \rightarrow \varphi))) = 1$.
 Let us prove that 

			\begin{equation}\label{eq:B2}
			Th(\mathcal{GK}^c),{^\ast\Box}^{=1}_u, \delta, \delta \rightarrow (\bigvee \Diamond^{<1}_u \rightarrow \varphi) \not \models_{[0,1]_G} \varphi
			\end{equation}			
			Suppose the contrary, with a view to contradiction.		
			Using completeness of $\models_{[0,1]_G}$ with respect to $\vdash_{\mathcal{G}}$, and applying Lemma \ref{lemma:propReduction} twice (once in each direction) and $M_\Box$ in between,  it follows that
			\begin{equation}\label{eq:contrad2}
			Th(\mathcal{GK}^c), \Box({^\ast\Box}^{=1}_u), \Box \delta, \Box (\delta \rightarrow ( \bigvee \Diamond^{<1}_u \rightarrow \varphi)) \models_{[0,1]_G} \Box \varphi
			\end{equation}
			But  this leads to a contradiction, since:
			\begin{itemize}
				\item $u(Th(\mathcal{GK}^c))= 1$ since $u \in W^\rho$ and $u(\Box({^\ast\Box}^{=1}_u))= 1$ by definition of ${^\ast\Box}^{=1}_u$. Moreover, $u(\Box \delta) > \alpha$ (Lemma \ref{lemma:deltaalpha}), and $u(\Box (\delta \rightarrow ( \bigvee \Diamond^{<1}_u \rightarrow \varphi))) = 1$ by assumption of the sub-case (B). Letting $\Gamma$ be the premises in  (\ref{eq:contrad2}), the previous amount to say that
				$\bigwedge_{\gamma \in \Gamma}  u(\gamma) > \alpha$. 
				
				\item However, $u(\Box \varphi) = \alpha$, contradicting \textit{3.} from Proposition \ref{prop:propStrongCompleteness}.
			\end{itemize}
			Thus, we have proven Condition (\ref{eq:B2}). This implies there exists an homomorphism $h \in Hom(\mathcal{L}_{\Box\Diamond}(V), [0,1]_G)$ that sends its premises to $1$ and the conclusion to some value strictly less than $1$. We claim this homomorphism $h$ meets the four conditions from the Proposition. Let us see why: 
			\begin{itemize}
				\item Since the premises in  (\ref{eq:B2}) are sent to $1$, we have that $h(Th(\mathcal{GK}^c)) = 1$ (proving \textit{(C1)}) and 
				$h({^\ast\Box}^{=1}_u) = 1$ (proving \textit{(C2)}),
				\item $h(\delta) = 1$ and $h(\varphi) < 1$, and from Observation \ref{obs:deltacond} these imply \textit{(C4)},
				
				\item $h(\delta \rightarrow(\bigvee \Diamond^{<1}_u \rightarrow \varphi )) = 1$, which together with $h(\delta) = 1$ and $h(\varphi) < 1$ imply $h(\bigvee \Diamond^{<1}_u) < 1$, namely, \textit{(C3)}. \qedhere
			\end{itemize}				
\end{proof}

It is easy that, since $u(\Box \top) = 1$, such an homomorphism further satisfies

\begin{itemize}[font=\itshape, align=left]
	\item[(C4.1)] $h(\varphi) < 1$.
\end{itemize}

The following are  some other properties of any $h$ as in the previous proposition. 
\begin{remark}\label{rem:properties}
	An homomorphism $h$ with properties \textit{(C1),(C2)} and \textit{(C3)} from Proposition \ref{prop:witnessBox}, further satisfies for any formulas $\theta_1, \theta_2, \theta \in \mathcal{L}_{\Box\Diamond}(V)$:
	\begin{itemize}[align=left]
		\item[\textit{(C2.a)}] $u(\Diamond \theta _{1})\leq u(\Box \theta_{2}) $ implies $h(\theta _{1})\leq h(\theta _{2})$ (since via $(FS2)$ $\theta_1 \rightarrow \theta_2 \in {^\ast\Box}^{=1}_u$);
		\item[\textit{(C2.b)}] For $\theta_1 \in \langle \rho \rangle$, $u(\Diamond \theta _{1}) < u(\Box \theta_{2}) $ implies $h(\theta_1) < h(\theta_2)$ (since $u(\Diamond \theta_1) < 1$ and $u(((\Box \theta_{2} ) \rightarrow \Diamond \theta _{1}) \rightarrow \Diamond \theta _{1}) = 1$ imply, via  $(T2)$, that $(\theta_2 \rightarrow \theta_1) \rightarrow \theta_1 \in {^\ast\Box}^{=1}_u$, and \textit{(C3)} further implies that $h(\theta_1) < 1$);
		\item[\textit{(C2.c)}] $0 < u(\Box \theta) $ implies $0 < h(\theta)$ (using \textit{(C2.b)}, since $u(\Diamond \bot) = 0$).
		\item[\textit{(C2.d)}]  $u(\Diamond \theta) = 0$ implies $0 = h(\theta)$ (using \textit{(C2.a)} since $u(\Diamond \theta) \leq h(\Box \bot)$).
	\end{itemize}
\end{remark}

Similarly to how it is done in \cite{CaRo15}, it is possible to build a G\"odel endomorphism that composed with the previous homomorphism will allow us to provide a world $v \in W^{\rho}$ such that $R^{\rho}uv$ and where $v(\varphi)$ is as near as possible to $u(\Box \varphi) = \alpha$. 

\begin{proposition}\label{prop:embeddingBox}
Let $\alpha < 1$, $\varphi  \in\Box^{=\alpha}_u$ and $\varepsilon > 0$. Then there is $w \in W^{\rho}$ such that $R^{\rho}uw$ and $w(\varphi) \in [\alpha, \alpha + \varepsilon]$.
\end{proposition}
\begin{proof}
	
	Let us consider the set $A =\{u(\Box \theta ):\theta \in \langle \rho \rangle\}$, and for any $a \in A$ let
	$h_a \coloneqq \bigwedge h(\Box^{=a}_u)$.\footnote{namely, $h_a = \min\{h(\theta) \colon \theta \in \langle \rho \rangle \text{ and } u(\Box \theta) = a\}$.}
Further let $h_\alpha^+ \coloneqq \min\{h_a \colon a \in A, h_\alpha < h_a\}$.\footnote{Recall that by convention the minimum of the empty set is $1$. We will be using this fact below without further notice.} From \textit{(C4.1)} we know $h_\alpha < h_\alpha^+$.
Observe condition \textit{(C4)} on $h$ implies directly that (*)$h_\alpha \leq h(\varphi) < h_a$ for any $\alpha < a$. Similarly, from \textit{(C2.c)} we get that (**)$h_0 < h_a$ for any $a > 0$.
	
	We inductively define the following sequence in $A$:
	\begin{align*}
	b_0 &\coloneqq \max \{a\in A \colon a<1\}, \\
	b_{i+1} & \coloneqq \max \{a\in A \colon a<b_{i} \text{ and }h_{a}<h_{b_{i}}\}.
	\end{align*}
	
	Since $A$ is finite, the previous is a strictly decreasing finite sequence, that has $0$ as last element (from observation (**) above).  Moreover, $\alpha = b_i$ for some $i$, from the observation (*) above. We will denote this index by $i_\alpha$ (i.e., $\alpha = b_{i_\alpha}$).
Also, by construction, the sequence $h_{b_{0}},h_{b_{1}},...$ is a strictly descending sequence with last element equal to $0$.

	Let us now define the sequence that will determine the upper bounds of our partial mappings. 
	\begin{align*}
	t_{0}&\coloneqq \bigwedge u(\Diamond \Diamond^{>b_0}_u), \\ 
	t_{i_{\alpha}} &\coloneqq (\alpha +\varepsilon )\wedge \bigwedge u(\Diamond \Diamond^{>\alpha}_u), \\
	t_{i+1} &\coloneqq b_{i} \wedge \bigwedge u(\Diamond \Diamond^{>b_{i+1}}_u) \text{ for } i+1 \neq i_\alpha
%
	\end{align*}
	By construction, $t_{i}>b_{i}$, which implies that also the sequence $\langle t_i\rangle$ is strictly decreasing.
	
	For simplicity in the notation, allow us to add the element $b_{-1} =1$ to the sequence (so $h_{b_{-1}} = 1$, from \textit{(C2)}). Then,
	Let $\sigma \colon [0,1] \rightarrow [0,1]$ be a strictly increasing uniform function with 
\begin{align*}	
\sigma([h_\alpha, h_\alpha^+)) &= [\alpha, t_{i_\alpha}), & \sigma([h_{b_i},h_{b_{i-1}})) &= [b_i,t_i) \text{ for } i \geq 0
	\end{align*}

	We know $t_{i}>b_{i}$ and $t_{i_\alpha} > \alpha$ by definition, so each interval in the right side is non-empty, proving $\sigma$ is well defined. Let us then denote $w \coloneqq \sigma \circ h$.
	Since $h_\alpha \leq h(\varphi) < h_\alpha^+$ (from \textit{(C4.1)}), and $t_{i_\alpha}\leq \alpha + \varepsilon$,  by definition we get $w(\varphi) \in [\alpha, \alpha + \varepsilon]$. We prove below that $R^{\rho}uw$, that is, for any formula $\theta \in \langle \rho \rangle$ we should see $u(\Box \theta) \leq w(\theta) \leq u(\Diamond \theta)$. This will conclude the proof of the Proposition.

	First, $w(\theta) = 1$ implies $h(\theta) = 1$. From \textit{(C3)} of $h$, this implies $u(\Diamond \theta) = 1$, and so, $w(\theta) \in [u(\Box \theta), u(\Diamond \theta)]$.
	
	If $w(\theta) <1$ then there is some $i \geq 0$ for which $h(\theta) \in [h_{b_{i}}, h_{b_{i-1}})$, since $h_{b_{-1}} = 1$ and the last element of the sequence is $0$. 
		\begin{itemize}
			\item $u(\Box \theta) \leq b_{i}$ follows from the definition of $b_{i}$. Indeed, $h_{u(\Box \theta)} \leq h(\theta) < h_{b_{i-1}}$  and $u(\Box \theta) < b_{i-1}$ (otherwise $b_{i-1} = u(\Box \theta)$, getting a contradiction with $h(\theta) < h_{b_{i-1}})$). Since $b_{i+1}$ is the maximum element in $A$ with those properties, $u(\Box \theta) \leq b_{i}$.
			\item To show the condition for $\Diamond$, we first prove that for any $\psi \in \langle \rho \rangle, j \geq 0$:
			\begin{equation*}
			u(\Diamond \psi) < t_{j} \text{ implies } h(\psi) \leq h_{b_{j}}
			\end{equation*}
			$u(\Diamond \psi) < t_{j}$ implies by definition of $t_j$ that $u(\Diamond \psi) \leq b_j$. From \textit{(C2.a)} we get $h(\psi) \leq h(\chi)$ for all $\chi  \in \langle \rho \rangle$ such that $u(\Box \chi) = b_j$. Thus, by definition, $h(\psi) \leq h_{b_j}$.
			
			\noindent Now, if $u(\Diamond \theta) < t_i$, from the previous equation we know that $h(\theta) \leq h_{b_{i}}$, so in fact, $h(\theta) = h_{b_{i}}$ and thus $w(\theta) = b_{i}$. But $h(\theta) = h_{b_{i}}$ implies $(1 >) b_i = u(\Box \theta) \leq u(\Diamond \theta)$ by $(T3)$, so $w(\theta) \leq u(\Diamond \theta)$
			
			 \noindent Otherwise, $t_{i} \leq u(\Diamond \theta)$.\qedhere			
		\end{itemize}
\end{proof}

The previous results will allow us to check the Truth-Lemma of the Canonical Model for formulas starting with $\Box$. We prove now results analogous to Propositions \ref{prop:witnessBox} and  \ref{prop:embeddingBox} but aiming towards the construction of a successor witnessing the values of the $\Diamond$ formulas.

Let us observe an easy fact on the behavior of formulas in $\Diamond^{<\alpha}_u$.

\begin{lemma}\label{lemma:deltaalphaDiamond}
	Let $\alpha >0$ and $\varphi  \in\Diamond^{=\alpha}_u$. Put
	\[\delta' \coloneqq (\varphi \to \bigvee \Diamond^{<\alpha}_u) \to \bigvee \Diamond^{<\alpha}_u\]
	Then $u(\Diamond \delta') =1$.
\end{lemma}
\begin{proof}
	It follows directly from $T^<_\Diamond$ and $K_\Diamond$, since $u(\Diamond \varphi) > u(\bigvee \Diamond \Diamond^{<\alpha}_u)$ by definition.
\end{proof}

\begin{proposition} \label{prop:witnessDiamond}
	Let $\alpha > 0$ and $\varphi \in \Diamond^{=\alpha}_u$. Then there exists $h \in Hom(\mathcal{L}_{\Box \Diamond}(V), [0,1]_G)$ such that
	\begin{enumerate}[align=left]
		\item[(C1)] $h(Th(\mathcal{GK}^c)) \subseteq\{1\}$,
		\item[(C2)] $h({^\ast\Box}^{=1}_u) \subseteq\{1\}$,
		\item[(C3)] $h(\psi) < 1$ for all $\psi \in \Diamond^{< 1}_u$,
		\item[(C4)$^\prime$] $h(\psi) < h(\varphi)$ for all $\psi \in \Diamond^{<\alpha}_u$.
	\end{enumerate}
\end{proposition}
\begin{proof}

Let us prove that 
\begin{equation}\label{eq:diamond}
Th(\mathcal{GK}^c),{^\ast\Box}^{=1}_u, \delta' \not \models_{[0,1]_G} \bigvee \Diamond^{< 1}_u
\end{equation}
Assume the contrary, with a view to contradiction. 
Similarly to the proofs from Proposition \ref{prop:witnessBox}, using completeness of $\models_{[0,1]_G}$,  Lemma \ref{lemma:propReduction}, the D.T for $\vdash_{\mathcal{GK}^c}$, the rule $M_\Box$ and theorems $K_\Diamond$ and $P$, and lastly again Lemma \ref{lemma:propReduction}, it follows that
\[Th(\mathcal{GK}^c),\Box {^\ast\Box}^{=1}_u  \models_{[0,1]_G} \Diamond \delta' \rightarrow \bigvee \Diamond \Diamond^{< 1}_u\]
This leads to a contradiction: on the one hand, $u(Th(\mathcal{GK}^c)) \subseteq\{1\}$ (because $u \in W^{\rho}$), $u(\Box {^\ast\Box}^{=1}_u) \subseteq\{1\}$ (by definition) and $u(\Diamond \delta') = 1$ from Lemma \ref{lemma:deltaalphaDiamond}; on the other hand, $u(\bigvee  \Diamond \Diamond^{< 1}_u) < 1$ by definition too, contradicting the definition of $\models_{[0,1]_G}$.

Condition \ref{eq:diamond} allows to conclude the proposition, since it implies there is 
$h \in Hom(\mathcal{L}_{\Box\Diamond}(V), [0,1]_G)$ evaluating the premises to $1$ and the conclusion to less than $1$. This $h$ satisfies the conditions of the Proposition, since:
\begin{itemize}
	\item Since the premises are sent to $1$, necessarily $h(Th(\mathcal{GK}^c))\subseteq\{1\}$ (namely,  \textit{(C1)})  and $h({^\ast\Box}^{=1}_u) \subseteq\{1\}$ (namely, \textit{(C2)}). 
	\item Since the conclusion is sent to less than $1$, we have that $h(\psi) < 1$ for all $\psi \in \Diamond^{< 1}_u$ (thus proving \textit{(C3)}),
	\item Using the last of the premises is sent to $1$ by $h$, i.e., $h(\delta') = 1$, it follows that $h((\varphi \to \bigvee \Diamond^{<\alpha}_u) \to \bigvee \Diamond^{<\alpha}_u) = 1$. From the previous point it follows that, in particular, $h(\bigvee \Diamond^{<\alpha}_u) < 1$. Thus necessarily $h(\bigvee \Diamond^{<\alpha}_u) < h(\varphi)$, proving \textit{(C4)$^\prime$}. \qedhere
\end{itemize}
\end{proof}

Homomorphism $h$ further satisfies (for formulas in $\langle \rho \rangle$) the conditions stated in Remark \ref{rem:properties}, since the necessary requisites are met. Moreover, in a dual way, and since $u(\Diamond \bot) = 0$, we now have that

\begin{itemize}[font=\itshape, align=left]
	\item[(C4.1)$^\prime$] $h(\varphi) > 0$.
\end{itemize}

We can again adapt the previous homomorphism in a dual way to how it was done in Proposition \ref{prop:embeddingBox}.
\begin{proposition}\label{prop:embeddingDiamond}
			Let $\alpha > 0$, $\varphi  \in\Diamond^{=\alpha}_u$ and $\varepsilon > 0$. Then there is $w \in W^{\rho}$ such that $R^{\rho}uw$ and $w(\varphi) \in [\alpha-\varepsilon, \alpha]$.
\end{proposition}
		\begin{proof}
		The proof is dual to the one of that Proposition  \ref{prop:embeddingBox}, swapping coherently $\Box$ and $\Diamond$ and handling the corresponding boundary values for $\Diamond$. We detail it for convenience of the reader.

			Let in this case $A =\{u(\Diamond \theta ):\theta \in \langle \rho \rangle\}$, and for any $a \in A$ let $h_a \coloneqq \bigvee h(\Diamond^{=a}_u)$.
		Further let $h_\alpha^- \coloneqq \max\{h_a \colon a \in A, h_a < h_\alpha\}$.\footnote{ Recall that by convention the maximum of the empty set is $0$. We will be using this fact below without further notice.} From \textit{(C4.1)$^\prime$} we know $h_\alpha > h_\alpha^-$.
		Observe condition \textit{(C4)$^\prime$} on $h$ further implies  that (*)$h_\alpha \geq h(\varphi) > h_a$ for any $\alpha > a$. Similarly, from \textit{(C3)} we get that (**)$h_a < h_1$ for any $a < 1$.
		
		We define the following sequence in $A$ (now, we start with the top boundaries, from below):
		\begin{align*}
		t_0 &\coloneqq \min \{a\in A \colon 0 < a\}, \\
		t_{i+1} & \coloneqq \min \{a\in A \colon a>t_{i} \text{ and }h_{a}>h_{t_{i}}\}.
		\end{align*}

		Since $A$ is finite, the previous is a strictly increasing finite sequence, that has $t_N=u(\Diamond \top)$ as last element.  Moreover, $\alpha = t_i$ for some $i$, from the observation (*). We will denote this index by $i_\alpha$ (i.e., $\alpha = t_{i_\alpha}$).

		Also, by construction, the sequence $h_{t_{0}},h_{t_{1}},\dots,h_{t_N}$ is a strictly increasing sequence with last element equal to $1$.

		Let us now define the sequence for the lower bounds.
		\begin{align*}
		b_{0}&\coloneqq \bigvee u(\Box \Box^{<t_0}_u), \\  
		b_{i_{\alpha}} &\coloneqq (\alpha -\varepsilon )\vee \bigvee u(\Box \Box^{<\alpha}_u), \\
		b_{i+1} &\coloneqq b_{i} \vee \bigvee u(\Box \Box^{<t_{i+1}}_u) \text{ for } i+1 \neq i_\alpha
		%
		\end{align*}
		
		By construction, $t_{i}>b_{i}$, which implies that also the sequence $\langle b_i\rangle$ is strictly increasing.
		
		For simplicity in the notation, allow us to add the element $t_{-1} =0$ to the sequence (so $h_{t_{-1}} = 0$ from \textit{(C2.d)}).
		Let then $\sigma \colon [0,1] \rightarrow [0,1]$ be a strictly increasing uniform function with 
		\begin{align*}	
		\sigma((h_\alpha^-, h_\alpha]) &= (b_{i_\alpha},\alpha], \\
		 \sigma((h_{t_i-1},h_{t_{i}}]) &= [b_i,t_i) \text{ for } i \geq 0
		\end{align*}

		We know $t_{i}>b_{i}$ and $b_{i_\alpha} < \alpha$ by definition, so each interval in the right side is non-empty, proving $\sigma$ is well defined. Let us then denote $w \coloneqq \sigma \circ h$.
		Since $h_\alpha \geq h(\varphi) > h_\alpha^-$ (from \textit{(C4)}$^{\prime}$), and $t_{i_\alpha}\geq \alpha - \varepsilon$,  by definition we get $w(\varphi) \in [\alpha-\varepsilon, \alpha]$. We prove below that $R^{\rho}uw$, that is, for any formula $\theta \in \langle \rho \rangle$ it holds $u(\Box \theta) \leq w(\theta) \leq u(\Diamond \theta)$. This concludes the proof of the Proposition.
		
		First, $w(\theta) = 0$ implies $h(\theta) = 0$. From \textit{(C2.c)} of $h$, this implies $u(\Box \theta) = 0$, and so, $w(\theta) \in [u(\Box \theta), u(\Diamond \theta)]$.
		
		If $w(\theta) >1$ then there is some $i \geq 0$ for which $h(\theta) \in (h_{t_{i-1}}, h_{t_{i}}]$, since $h_{t_{-1}} = 0$ and the last element of the sequence is $1$. 
		\begin{itemize}
			\item $u(\Diamond \theta) \geq t_i$ follows from the definition of $t_i$. Indeed, $h_{u(\Diamond \theta)} \geq h(\theta) > h_{t_{i-1}}$, and $u(\Diamond \theta) > t_{i-1}$ (otherwise $t_{i-1} = u(\Diamond \theta)$, getting a contradiction with $h(\theta) > h(t_{b_{i-1}})$). Since $t_{i}$ is the minimum element in $A$ with those properties, $u(\Diamond \theta) \geq t_{i}$.
			\item To show the condition for $\Box$, we first prove that for any $\psi \in \langle \rho \rangle, j \geq 0$:
			\begin{equation*}
			u(\Box \psi) > b_{j} \text{ implies } h(\psi) \geq h_{t_{j}}
			\end{equation*}
			$u(\Box \psi) > b_{j}$ implies by definition of $b_j$ that $u(\Box \psi) \geq t_j$. From \textit{(C2.a)} we get $h(\psi) \geq h(\chi)$ for all $\chi  \in \langle \rho \rangle$ such that $u(\Diamond \chi) = t_j$. Thus, by definition, $h(\psi) \geq h_{t_j}$.

			\noindent Now, if $u(\Box \theta) > b_i$, from the previous equation we know that $h(\theta) \geq h_{t_{i}}$, so in fact, $h(\theta) = h_{t_{i}}$ and thus $w(\theta) = t_{i}$. But $h(\theta) = h_{t_{i}}$ implies $(0 <) t_i = u(\Diamond \theta) \geq u(\Box \theta)$ by $(T3)$, so $w(\theta) \geq u(\Box \theta)$
			
			\noindent Otherwise, $b_{i} \geq u(\Box \theta)$.
			
		\end{itemize}
		
\end{proof}

		With the previous machinery, we can now go back to prove that the model $\m{M}^{\rho}$ is indeed canonical for formulas in $\langle \rho \rangle$.
	
\begin{lemma} [Truth-Lemma] \label{lemma:truthLemma} $e^{\rho}(u,\varphi )=u(\varphi )$ for any $\varphi \in \langle \rho \rangle$ and any $u\in W^{\rho}$.
\end{lemma}

\begin{proof}
	This can be proven, as usual, by induction on the complexity of the formulas. Propositional cases are trivial, and so the relevant cases are the steps of modal operations. Thus, applying Induction Hypothesis, the objective is to prove that for any $\Box \varphi, \Diamond \varphi \in \langle \rho \rangle$:
	\[ (TL_\Box): u(\Box \varphi ) = \bigwedge_{R^\rho uv} v(\varphi)\quad  \text{ and } \quad (TL_\Diamond): u(\Diamond \varphi )   = \bigvee_{R^\rho uv} v(\varphi)\]

	$\leq$ in $(TL_\Box)$ and $\geq$ in $(TL_\Diamond)$ follow immediately by definition.
	These inequalities further proof the full equality whenever $u(\Box \varphi ) = 1$ or $u(\Diamond \varphi ) = 0$, respectively.
	
	For the rest of the cases,  Proposition \ref{prop:embeddingBox} proves that $\bigwedge_{R^\rho uv} v(\varphi) \leq u(\Box \varphi)$, and correspondingly, Proposition \ref{prop:embeddingDiamond} proves $u(\Diamond \varphi ) \leq \bigvee_{R^\rho uv} v(\varphi)$.
\end{proof}

\begin{theorem} [Weak completeness]
\label{JointCompleteness} For any formula $\varphi $ in $\mathcal{L}_{\Box \Diamond }$%
 \[\vdash _{\mathcal{GK}^c}\varphi \text{ if and only if } \models_{\m{GK}^c} \varphi \]
\end{theorem}

\begin{proof}
	Soundness of the axioms is simple to check. For what concerns completeness,
	 assume  $\not\vdash _{\mathcal{GK}^c}\varphi .$ Then $Th(\mathcal{GK}^c)\not\vdash _{\mathcal{G}} \varphi $ by Lemma \ref{lemma:propReduction}, and thus there is, by Proposition \ref{prop:propStrongCompleteness}, an element $v \in Hom(\mathcal{L}_{\Box\Diamond}(V), [0,1]_G)$ such that $v(\varphi )< v(Th(\mathcal{GK}^c)) = \{1\}.$ Then $v$ is a world of
the canonical model $\m{M}^\varphi$ and by Lemma \ref{lemma:truthLemma}, $e^{\varphi}(v,\varphi )=v(\varphi )<1.$ Thus $\not\models_{\m{GK}^c} \varphi .$
\end{proof}

The previous proof of completeness for theorems of the logic, together with the Deduction Theorem for $\mathcal{GK}^c$ allows us to easily generalize the completeness result to deductions.

\begin{corollary} [Finite strong completeness]
\label{JointCompleteness2} For any finite set of formulas $\Gamma \cup \{\varphi\} \subseteq \mathcal{L}_{\Box\Diamond}$ the following are equivalent:
 \[\Gamma \vdash _{\mathcal{GK}^c}\varphi \text{ if and only if } \Gamma \models_{\m{GK}^c} \varphi \]
\end{corollary}
\begin{proof}
	
	Left-to-right is a direct consequence of the D.T. For the other direction, assume $\Gamma \not \vdash_{\mathcal{GK}^c} \varphi$. Thus, from Lemma \ref{lemma:propReduction} and completeness of G\"odel propositional logic (Prop. \ref{prop:propStrongCompleteness}) there is $h \in Hom(\mathcal{L}_{\Box\Diamond}, [0,1]_G)$ such that
	 $h(\Gamma) \subseteq\{ 1\}$, $h(Th(\mathcal{GK}^c)) \subseteq \{1\}$ and $h(\varphi) < 1$. In particular, $h$ is an element of the universe of the  canonical model of $\bigwedge \Gamma \rightarrow \varphi$ (or equivalently, of any formula containing both $\Gamma$ and $\varphi$ in its set of subformulas). From the Truth Lemma \ref{lemma:truthLemma} we know that in this model $e(h, \bigwedge\Gamma) = h(\bigwedge\Gamma) = 1$, and $e(h, \varphi)= h(\varphi) < 1$, proving that $\Gamma \not \models_{\m{GK}^c} \varphi$.
\end{proof}

Moreover, as it is done in \cite[Theorem~3.1]{CaRo15}, it is possible to extend this completeness  to infinite sets of formulas, as long as they are built on a countable set of variables.

\begin{corollary} [Strong completeness]
For any countable set of formulas $\Gamma$ and formula $\varphi \in \mathcal{L}_{\Box\Diamond}(V)$,
 \[\Gamma \vdash _{\mathcal{GK}^c}\varphi \text{ if and only if } \Gamma \models_{\m{GK}^c} \varphi \]
\end{corollary}
\begin{proof}
The proof is almost the same as in \cite[Theorem~3.1]{CaRo15}, only taking into account that, when building the theory that models the class $\m{GK}^c$ inside classical first order logic, we need to restrict the value of the accessibility relation to $\{0, 1\}$. We briefly reproduce the proof here for convenience of the reader. Some familiarity with first order logic is assumed.

Let $\Gamma$ be countable and $\Gamma \nvdash_{\mathcal{GK}^c}\varphi $ \ and consider the first order theory $\Gamma^{\ast }$ with two unary relation symbols\ $W,P,$ a binary relation symbol $<$, three constant symbols $0,1,c,$ two binary function symbols $\circ ,S,$ and a unary function symbol $f_{\theta }$ for each $\theta \in \mathcal{L}_{\square \Diamond }(V),$ where $V$ is the set of propositional variables
occurring in formulas of $\Gamma,$ and having for axioms:

\medskip

$\forall x\lnot (W(x) \wedge P(x))$

$\forall x (W(x) \vee \lnot W(x))$

``$(P,<)$ is a strict linear order with minimum $0$ and maximum $1"$

$\forall x\forall y(W(x)\wedge W(y)\rightarrow (S(x,y)= 1 \vee S(x,y) = 0))$

$\forall x\forall y(P(x)\wedge P(y)\rightarrow (x\leq y\wedge x\circ
y=1)\vee (x>y\wedge x\circ y=y))$

$\forall x(W(x)\rightarrow f_{\bot }(x)=0)$

for each $\theta ,\psi \in \mathcal{L}_{\square \Diamond }$ the sentences:

$\forall x(W(x)\rightarrow P(f_{\theta }(x)))$

$\forall x(W(x)\rightarrow f_{\theta \wedge \psi }(x)=\min \{f_{\theta
}(x),f_{\psi }(x)\})$

$\forall x(W(x)\rightarrow f_{\theta \rightarrow \psi }(x)=(f_{\theta
}(x)\circ f_{\psi }(x))$

$\forall x(W(x)\rightarrow f_{\square \theta }(x)=\inf_{y}(S(x,y)\circ
f_{\theta }(y))$

$\forall x(W(x)\rightarrow f_{\Diamond \theta }(x)=\sup_{y}(\min
\{S(x,y),f_{\theta }(y)\})$

for each $\gamma \in \Gamma$ the sentence: $f_{\gamma }(c)=1$

finally, $W(c)\wedge (f_{\varphi }(c)<1).$

\medskip

\noindent For each finite part $t$ of $\Gamma^{\ast }$ let $F$ be a finite fragment of $\mathcal{L}_{\square \Diamond }$ containing $\{\theta
:f_{\theta }\ $occurs in $t\}.$ Since $F\cap \Gamma \nvdash_{\mathcal{GK}^c}\varphi $ by hypothesis, then, by weak completeness, there is a crisp GK-model $\,M_{F}=(W,S^{F},e^{F})$ and $c\in W$ such that $e^{F}(c,\theta )=1$ for each $\theta \in F\cap \Gamma$ and $e^{F}(c,\varphi )<1.$ Therefore the first order structure $(W\sqcup \lbrack 0,1],W,[0,1],<,0,1,c,\Rightarrow ,S^{F},\{f_{\theta }\}_{\theta \in \mathcal{L}_{\square \Diamond }}),$ with $f_{\theta }:W\rightarrow \lbrack 0,1]$ defined as $f_{\theta }(x)=e^{F}(x,\theta ),$ is clearly a model of $t.$ By
compactness of first order logic and the downward Löwenheim theorem, $\Gamma^{\ast }$ has a countable model $M^{\ast }=(B,W,P,<,0,1,c,\circ
,S,\{f_{\theta }\}_{\theta \in \mathcal{L}_{\square \Diamond }}).$ Using Horn's lemma \cite{Horn69}, $(P,<)$ may be embedded in $(\mathbb{Q}\cap \lbrack
0,1],<)$ preserving $0,1,$ and all suprema and infima existing in $P$; therefore, we may assume without loss of generality that the
function $S$ is crisp and the ranges of the $f_{\theta }'s$ are contained in $[0,1].$ Then, it is straightforward to verify that $M=(W,S,e),$ where $e(w,\theta )=f_{\theta }(w)$ for all $w\in W$ and\ $\theta \in \mathcal{L}_{\square \Diamond }(V),$ is a crisp GK-model with a distinguished world $c$ such that $M, c \models \Gamma,$ and $M, c \not\models \varphi .$ Hence, $\Gamma\not\models _{{\m{GK}^c}}\varphi .$ 
\end{proof}

Strong completeness allows us to easily prove that extending $\mathcal{GK}^c$ with the unrestricted necessity rule ${^\ast N}_\Box\colon  \varphi \vdash \Box \varphi$ (i.e., affecting all formulas and not only theorems of the logic) provides a complete axiomatization of $\models^g_{\m{GK}^c}$.
Let us denote this axiomatic system by $^\ast \mathcal{GK}^c$. Moreover, for any formula $\psi$ we will write
\begin{align}
\Box^0 \psi \coloneqq& \psi, & \Box^{k+1}\psi \coloneqq \Box \Box^k\psi
\end{align}
and the corresponding analogous meaning for what concerns sets of formulas.

\begin{lemma}
	For any set $\Gamma \cup \{\varphi\}\subseteq \mathcal{L}_{\Box\Diamond}$
	 \[\Gamma \vdash _{^\ast \mathcal{GK}^c}\varphi \text{ if and only if } \{\Box^k \Gamma \}_{k \in \omega} \vdash_{\mathcal{GK}^c} \varphi \]
\end{lemma}
\begin{proof}
	Right-to-left direction is immediate. For what concerns the other we can simply reason by induction on the length of the derivation of $\varphi$ from $\Gamma$. The step for M.P. is immediate. For the ${^\ast N}_\Box$ step, assume $\Gamma \vdash _{^\ast \mathcal{GK}^c}\psi$. By Induction Hypothesis, we know $\{\Box^k \Gamma \}_{k \in \omega} \vdash_{\mathcal{GK}^c}  \psi$. Now, applying $(M_\Box)$ we get
	 $\Box \{\Box^k \Gamma \}_{k \in \omega} \vdash_{\mathcal{GK}^c} \Box \psi$. Since $\Box \{\Box^k \Gamma \}_{k \in \omega} \subset \{\Box^k \Gamma \}_{k \in \omega}$, this concludes the proof.
\end{proof}

\begin{corollary}[Strong global completeness]
	For any set of formulas $\Gamma \cup \{\varphi\} \subseteq \mathcal{L}_{\Box\Diamond}$,
	\[\Gamma \vdash _{^\ast \mathcal{GK}^c}\varphi \text{ if and only if } \Gamma \models_{\m{GK}^c}^g \varphi \]
\end{corollary}
\section{Some axiomatic extensions}\label{sec:extensions}
In a similar fashion to how it is done in \cite{CaRo15}, it is easy to axiomatize some of the better known frame structural properties. Moreover, since the accessibility relation in the models from $\m{GK}^c$ is classical, it is possible to also address some properties whose characterization in the full $\mathcal{GK}$ is unknown (eg. \textit{seriality}, $\forall x \exists y R(x,y)$).

The canonical model built in the previous sections is, as it happens in \cite{CaRo15}, determined by a finite subset of formulas, and in that sense, it lacks optimality with respect to the accessibility relation. 
Namely, its accessibility relation can be further extended in such a way that the Truth Lemma (Lemma \ref{lemma:truthLemma}) holds for all formulas in the language (not only those in $\langle \rho \rangle$). The procedure in order to do so is very similar to the one in the above reference, only taking into account the restriction to crisp models.

\begin{definition}[\textit{c.f.} Def.~4.1 from \cite{CaRo15}]
Given a crisp G\"odel-Kripke model $\m{M} = \langle W, R, e\rangle$ we define $R^+ \coloneqq \{\langle v,w\rangle \colon e(v, \Box \varphi) \leq e(w,\varphi) \leq e(w, \Diamond \varphi) \text{ for all }\varphi \in  \mathcal{L}_{\Box\Diamond}\}$, and denote $\m{M}^+ \coloneqq \langle W, R^+, e\rangle$. We call $\m{M}$ \termDef{optimal} if $R^+ = R$.
\end{definition}

For simplicity in the notation, we will write $e^+$ when evaluating in the extended model.
 It is easy to see by induction on the formulas that any model is equivalent to an optimal one.
\begin{lemma}\label{lem:optimalModel}
For any model $\m{M} \in \m{GK}^c$ and any formula $\varphi \in  \mathcal{L}_{\Box\Diamond}$,
\[e(v,\varphi) = e^+(v, \varphi).\]
\end{lemma}
 \begin{proof}
 Clearly $R \subseteq R^+$, so $e(v,\Box \varphi) \geq e^+(v, \Box \varphi)$ for every formula $\varphi$.
 On the other hand, if $R^+(v,w)$, by definition $e(v,\Box \varphi) \leq e(w, \varphi)$.
 Thus, $e(v, \Box \varphi) \leq \bigwedge_{R^+(v,w)} e(w, \varphi) = e^+(v,\Box \varphi)$. The proof for $\Diamond$ formulas is analogous.
 \end{proof}
The previous result provides completeness of $\vdash_{\mathcal{GK}^c}$ with respect to optimal models. 
One of their benefits is that it is easier to check the correspondence between axiomatic extensions and frame conditions over them.

Usual definitions of the classical frame conditions are preserved for crisp G\"odel-Kripke models, given that the frame itself is classical. Thus, a frame (and correspondingly, a model) is \textit{reflexive, transitive, symmetric, euclidean} and \textit{serial} if the corresponding classical conditions hold over its accessibility relation.

The modal schemes that characterize classically the previous conditions are the following:\footnote{In the classical case, only one of the two modal schemes corresponding to each condition is necessary to characterize the class of frames, due to the inter-definability of the modalities.}
\begin{align*}
T_\Box &\  \Box \varphi \rightarrow \varphi & T_\Diamond &\  \varphi \rightarrow \Diamond \varphi & ((\mathcal{M})\ reflexivity)\\
4_\Box &\  \Box \varphi \rightarrow \Box \Box \varphi & 4_\Diamond &\  \Diamond \Diamond \varphi \rightarrow \Diamond \varphi & ((\emph{4})\ transitivity)\\
B_1 &\  \varphi \rightarrow \Box \Diamond \varphi & B_2 &\  \Diamond \Box \varphi \rightarrow \varphi & ( (\mathcal{B})\ symmetry)\\
5_1 &\  \Diamond \varphi \rightarrow \Box \Diamond \varphi & 5_2 &\  \Diamond \Box \varphi \rightarrow \Box \varphi & ((\emph{5})\ euclideanicity)\\
D &\  \Diamond \top & & & ((\mathcal{D})\ seriality)
\end{align*}

Validity of the previous pairs of axioms in the corresponding classes of models is direct. It is also not hard to see that, over optimal models, the axiom schemata indeed characterize the corresponding frame conditions.
\begin{lemma}\label{lem:frameCorrespondance}
Let $\m{M}$ be an optimal crisp G\"odel-Kripke model. Then the following hold:
\begin{itemize}
\item $\m{M}$ is \emph{reflexive} if and only if it validates the schemes $T_\Box, T_\Diamond$;
\item $\m{M}$ is \emph{transitive} if and only if it validates the schemes $4\Box, 4_\Diamond$;
\item $\m{M}$ is \emph{symmetric} if and only if it validates the schemes $B_1, B_2$;
\item $\m{M}$ is \emph{euclidean} if and only if it validates the schemes $5_1, 5_2$;
\item $\m{M}$ is \emph{serial} if and only if it validates the scheme $D$;
\end{itemize}
\end{lemma}
\begin{proof}
The proof is a simple application of the definition of optimal model. Readers interested in the details can consult of \cite[Prop.~4.1]{CaRo15}. The case of seriality follows by definition (and holds also for non-optimal models).
\end{proof}

Combining the completeness result with respect to optimal models pointed out in Lemma \ref{lem:optimalModel} and the previous lemma, the axiomatization of the analogous extensions of $\mathcal{GK}^c$ is clear.

\begin{theorem}
Let $P \subseteq \{\{T_\Box, T_\Diamond\}, \{4_\Box, 4_\Diamond\},  \{B_1, B_2\}, \{5_1, 5_2\}, D\}$, and $\mathcal{P}$ be the corresponding set of frame restrictions associated to the axioms from $P$ (in the sense of \ref{lem:frameCorrespondance}). Consider $\mathcal{GK}^c_P$ be the axiomatic system $\mathcal{GK}^c$ extended by the axiom schemata from $\bigcup P$. Then $\mathcal{GK}^c_P$ is strongly complete with respect to the G\"odel Kripke models with frame conditions from $\mathcal{P}$.
\end{theorem}

We would like to close this section pointing out that the previous result implies that $\mathcal{GK}^c_P$ with $P = \{\{T_\Box, T_\Diamond\}, \{4_\Box, 4_\Diamond\}, \{5_1, 5_2\}\}$ is an axiomatization of the monadic G\"odel logic $S5(\mathcal{G})$, namely, the one arising from G\"odel-Kripke models where the accessibility relation is a equivalence relation. This logic was studied by H\'ajek in \cite{Ha2010} but his proof of completeness is wrong, as it is pointed out in \cite{Casetal2020}, where another alternative proof is provided. In addition, it is worth mentioning that $\mathcal{GK}^c_P$ with $P$ as above is equivalent to\footnote{We thank Xavier Caicedo for this observation.} the logic $GS5 + C$ below (studied in \cite{CaRo15}):
\[
(GS5) \left \{
\begin{array} {l l}
\mathcal{G} & \text{Propositional G\"odel logic} \\
K_\Box &  \Box(\varphi \to \psi) \to \Box \varphi \to \Box \psi \\
K_\Diamond & \Diamond(\varphi \vee \psi) \to \Diamond \varphi \vee \Diamond \psi \\
P & \Box(\varphi \to \psi) \to \Diamond \varphi \to \Diamond \psi \\
T & \Box \varphi \to \varphi \text{ and } \varphi \to \Diamond \varphi \\
5 & \Diamond \varphi \to \Box \Diamond \varphi \text{ and } \Diamond \Box \varphi \to \Box \varphi \\
\end{array}
\right \}
 + 
(C)\ \Box(\Box \varphi \vee \psi) \to \Box \varphi \vee \Box \psi \\
\]

\section{Non Interdefinability of the modal operators}\label{sec:nonInterdef}

One might wonder if it is possible to define $\Box$ from $\Diamond$ or vice-versa in the logic arising from the class of crisp models $\m{GK}^c$. While it is easy to see that the usual definition of one modality from the other using negation as it is done in classical modal logic does not hold (see eg. \cite{CaRo10}), it could be the case that other  possible formulations did (for instance, in \cite{ViEsGo19} it is proven how this can be done in  logics with canonical constants). While the failure for inter-definability of $\Box$ and $\Diamond$ was expected, up to our knowledge there were no proofs in the literature showing this was indeed the case.

In this section, we  prove that indeed it is not
possible to define $\Box$ from $\Diamond$ or vice-versa, showing that the axiomatization $\mathcal{GK}^c$ we provided in Section \ref{sec:axiomatization}  is a new logic different from $\mathcal{GK}_\Box^c$ and  $\mathcal{GK}_\Diamond^c$ the fragments. This implies that the modalities are not interdefinable neither in the larger class of models $\m{GK}$, where the accessibility relation is $[0,1]$-valued.

First observe it is clear that $\Diamond$ cannot be possibly defined from $\Box$ in $\m{GK}$, because the $\Box$-fragment over $\m{GK}$ is complete with respect to $\m{GK}^c$, while this is not the case for the $\Diamond$-fragment \cite{CaRo10}. The case of $\Box$ not being definable from $\Diamond$ in $\m{GK}$ will follow from the same result over $\m{GK}^c$, which we prove below.

 To be precise, we say that $\Box$ is definable from $\Diamond$ in the class of models $\m{C}$ if there is some $\Box$-free formula
 $\phi_\Box(x)$ such that,
 for any  $\Box$-free formula $\varphi$, any model $\m{M} \in \m{C}$ and any $v \in W$,
 \[e(v, \Box \varphi) = e(v, \phi_\Box(\varphi))\]
Dually, $\Diamond$ is definable from $\Box$ if there is some $\Diamond$-free formula $\phi_\Diamond$ such that  for any  $\Diamond$-free formula $\varphi$, any model $\m{M} \in \m{C}$ and any $v \in W$,
 \[e(v, \Diamond \varphi) = e(v, \phi_\Diamond(\varphi))\]

 For $\m{C} = \m{GK}^c$, via the completeness result proven before (Theorem \ref{JointCompleteness}), the previous definitions are equivalent to say that there is a $\Box$-free formula $\phi_\Box(x)$ such that
$\Box \varphi \leftrightarrow \phi_\Box(\varphi)$ is a theorem of the logic $\mathcal{GK}^c$, for any $\Box$-free formula $\varphi$ (and the dual for the definition of $\Diamond$ formulas in terms of formulas with only $\Box$).

 \begin{lemma}
$\Box$ is not definable from $\Diamond$ in $\m{GK}^c$, and $\Diamond$ is not definable from $\Box$ in $\m{GK}^c$.
 \end{lemma}
\begin{proof}

\

\noindent
\begin{minipage}{0.64\textwidth}

In order to do so, we will define a $\mathcal{GK}^c$ algebra $\alg{A}$, and choose, for each one of the claims  above, a corresponding reduct $\alg{A}_\Diamond$ (correspondingly $\alg{A}_\Box$) such that $\alg{A}_\Diamond$ is a $\mathcal{GK}^c_\Diamond$ algebra but not a $\mathcal{GK}^c$-subalgebra of $\alg{A}$ (correspondingly, an $\mathcal{GK}^c_\Box$ algebra that is  not a $\mathcal{GK}^c$-subalgebra). These prove the lemma: if $\Box$ could be defined from $\Diamond$ in the above sense, since the formula $\Box\varphi \leftrightarrow \phi_\Box(\varphi)$ should be valid in $\alg{A}$, necessarily also $\alg{A}_\Diamond$ should be a $\mathcal{GK}^c$ algebra (and dually for $\Diamond$).

In order to build the above algebra, consider the frame $\m{F}$ in the right side.
 \end{minipage}
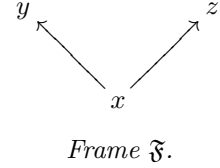
\begin{minipage}{0.33\textwidth}
\[
\xymatrix{%
   y & & z \\
    & x \ar[ul] \ar[ur] & }
\]
\center{\emph{Frame $\m{F}$.}}
\end{minipage}
The complex algebra arising from $\m{F}$ (see eg. \cite{CaRo15} for the general construction) is the algebra\footnote{In order to lighten the notation, we will denote a function $g \in [0,1]^{\{x,y,z\}}$ simply by the tripla $\langle g(x){,} g(y){,} g(z) \rangle \in [0,1]^3$.}
\[\alg{A} \coloneqq \langle [0,1]^3, \wedge, \vee, \rightarrow, \Box, \Diamond, 0, 1\rangle\]
where $\wedge, \vee, \rightarrow$ are interpreted as the G\"odel operations component-wise, $0$ and $1$ are the constant mappings to $0$ and $1$ respectively, and
\[\Box \langle a{,}b{,}c \rangle \coloneqq \langle b \wedge c{,} 1{,} 1 \rangle, \qquad  \Diamond \langle a{,}b{,}c \rangle \coloneqq \langle b \vee c{,} 0{,} 0 \rangle
\]

To prove that $\Box$ is not definable from $\Diamond$, let $\alg{A}_\Diamond$ be the $\Diamond$-subalgebra of $\alg{A}$ generated by the element $\langle 0{,} \frac{1}{2}{,} \frac{1}{3} \rangle$. It can be checked that the universe of $\alg{A}_\Diamond$ is the set
\[\{\langle 0{,}0{,}0\rangle{,} \langle 1{,}1{,}1\rangle{,} \langle 1{,}0{,}0\rangle{,}\langle 0{,}1{,}1\rangle{,} \langle 0{,}\tfrac{1}{2}{,}\tfrac{1}{3}\rangle{,} \langle 1{,}\tfrac{1}{2}{,}\tfrac{1}{3}\rangle{,} \langle \tfrac{1}{2}{,}\tfrac{1}{2}{,}\tfrac{1}{3}\rangle{,}\langle \tfrac{1}{2}{,}0{,}0\rangle{,} \langle \tfrac{1}{2}{,}1{,}1\rangle\}\]

This is a G\"odel subalgebra of $\alg{A}$ that is further closed under applications of $\Diamond$, thus a $\mathcal{GK}^c_\Diamond$-algebra.
However,  $\Box \langle 0{,} \frac{1}{2}{,} \frac{1}{3} \rangle = \langle \frac{1}{3}{,} 1{,} 1\rangle$ is not an element of $\alg{A}_\Diamond$, so it is not a $\mathcal{GK}^c$ subalgebra of $\alg{A}$.

To prove that $\Diamond$ is not definable from $\Box$, let $\alg{A}_\Box$ be the $\Box$-subalgebra of $\alg{A}$ generated by the same element from before, $\langle 0, \frac{1}{2}, \frac{1}{3} \rangle$. It is a matter of calculations to see that the universe of $\alg{A}_\Box$ is the set \[\{\langle 0{,}0{,}0\rangle, \langle 1{,}1{,}1\rangle{,} \langle 1{,}0{,}0\rangle{,}\langle 0{,}1{,}1\rangle{,}\langle 0{,}\tfrac{1}{2}{,}\tfrac{1}{3}\rangle{,} \langle 1{,}\tfrac{1}{2}{,}\tfrac{1}{3}\rangle{,} \langle \tfrac{1}{3}{,}\tfrac{1}{2}{,}\tfrac{1}{3}\rangle{,} \langle \tfrac{1}{3}{,}0{,}0\rangle{,}\langle \tfrac{1}{3}{,}1{,}1\rangle\}\]
This is a G\"odel subalgebra of $\alg{A}$ that is further closed under applications of $\Box$, thus a $\mathcal{GK}^c_\Box$-algebra.

On the other hand, $\Diamond \langle 0{,} \frac{1}{2}{,} \frac{1}{3} \rangle = \langle \frac{1}{2}{,} 0{,} 0\rangle$ is not an element in the previous $\alg{A}_\Box$, proving that in $\alg{A}_\Box$  is not a $\mathcal{GK}^c$ subalgebra of $\alg{A}$.

\end{proof}

Since $\m{GK}^c \subset \m{GK}$, the following is immediate.

 \begin{corollary}
$\Box$ is not definable from $\Diamond$ in $\m{GK}$.
 \end{corollary}

\section{Conclusions}
In this work, we have studied G\"odel many-valued logics extended with modal operators $\Box$ and $\Diamond$
interpreted over the class of models $\m{GK}^c$. The main contribution of this paper has been establishing an axiomatization  $\mathcal{GK}^c$ strongly complete with respect to the logic of $GK^c$-models. The proposed axiomatization is built by taking the one introduced in \cite{CaRo15} plus a simple axiom coming from the study of positive modal logics. In addition, we prove that both $\Box$-fragment and $\Diamond$-fragment are strictly included in our logic $\mathcal{GK}^c$.\\

We leave some open questions concerning the studied framework:
\begin{enumerate}
\item Is the axiom $Cr$ derivable from the axioms of $\mathcal{GK}$ and the rule $R_\Diamond$?
\item We know our logic $\mathcal{GK}^c$ is decidable under the alternative semantics proposed in \cite{CaMe13}, but is it possible to extend this result to other logics whose accessibility relations satisfy conditions such as reflexivity, symmetry, and transitivity?
\item What is the computational complexity of validity in $\m{GK}^c$?
\end{enumerate}

\noindent
\textbf{Acknowledgments}
The authors are thankful to the anonymous reviewers for their useful comments,  that have helped to improve the layout of the paper. 
This project has received funding from the following sources: 1) the European Union's Horizon 2020 Research and Innovation program under the Marie Sklodowska-Curie grant agreement No 689176 (SYSMICS project); 2) the grant no. CZ.02.2.69/0.0/0.0/17\_050/0008361 of the Operational programme Research, Development, Education of the Ministry of Education, Youth and Sport of the Czech Republic, co-financed by the European Union; 3) the Spanish MINECO project RASO (TIN2015-71799-C2-1-P) and 4) the Argentinean project PIP CONICET 11220150100412CO and UBA-CyT 20020150100002BA.

%
%
%
\end{document}